\documentclass[a4paper,12pt]{article}
\usepackage{fancyhdr,graphicx}
\usepackage{amsfonts}
\usepackage{amssymb}
\usepackage{amsthm}
\usepackage{newlfont}
\usepackage{amsmath}
\usepackage{mathrsfs}
\usepackage[top=2.5cm,bottom=3.0cm,left=2.5cm,right=2.5cm]{geometry}
\usepackage{multirow}
\usepackage{array}
\usepackage{longtable}
\usepackage{bm}
\usepackage{setspace}

\newcommand{\trou}{\vspace{1 mm}}
\newcommand{\noi}{\noindent}

\usepackage{latexsym,bm}
\title{\vspace{0cm}{\Huge {Proofs of two conjectures on generalized Fibonacci cubes\thanks{This
work is supported by NSFC (Grant no. 61073046).}}}}

\author{Jianxin Wei$^{a,b}$, Heping Zhang $^{a,}$\footnote{Corresponding author.} \\
\footnotesize {$^{a}$ School of Mathematics and Statistics, Lanzhou University, Lanzhou, Gansu 730000, P. R. China}\\
\footnotesize{$^{b}$School of Mathematics and Statistics Science, Ludong University, Yantai, Shandong 264025, P. R. China}
\\\small{E-mail addresses: wjx0426@163.com, zhanghp@lzu.edu.cn}}

\date{}

\begin{document}

\maketitle

\begin{abstract}
A binary string $f$ is a factor of string $u$ if $f$ appears as a sequence of $|f|$ consecutive bits of $u$, where $|f|$ denotes the length of $f$.
Generalized Fibonacci cube $Q_{d}(f)$ is the graph obtained from the $d$-cube $Q_{d}$ by removing all vertices that contain a given binary string $f$ as a factor.
A binary string $f$ is called good if $Q_{d}(f)$ is an isometric subgraph of  $Q_{d}$ for all $d\geq1$, it is called bad otherwise.
The index of a binary string $f$, denoted by $B(f)$,
is the smallest integer $d$ such that $Q_{d}(f)$ is not an isometric subgraph of $Q_{d}$.
Ili\'{c}, Klav\v{z}ar and Rho conjectured that $B(f)<2|f|$ for any bad string $f$.
They also conjectured that if $Q_{d}(f)$ is an isometric subgraph of $Q_{d}$,
then $Q_{d}(ff)$ is an isometric subgraph of $Q_{d}$.
We confirm the two conjectures by obtaining a basic result:
if there exist $p$-critical words for $Q_{B(f)}(f)$, then $p$=2 or $p=3$.

\end{abstract}

\section{Introduction}

Fibonacci cube $\Gamma_{d}$ was introduced by Hsu \cite{whs} as a model for interconnection networks,
which has beneficial properties similar to those of the hypercube.
It can be seen as the subgraph obtained from hypercube $Q_{d}$ by removing
all the vertices that contain two consecutive 1s.
As Klav\v{z}ar and \v{Z}igert \cite{sp} showed that Fibonaccenes have resonance graphs that are exactly the Fibonacci cubes.
More generally, Zhang et al. \cite{zoy} described the class of
planar bipartite graphs that have Fibonacci cubes as their resonance graphs.
Taranenko and Vesel \cite{tara} showed that
Fibonacci cubes can be recognized in $O(m$log$n)$
time (where $m$ is the size and $n$ the order of a given graph).
There are more results on the structural properties of  Fibonacci cubes, see \cite{dt, Gre, sk1, ms},
and \cite{sk} for a recent survey.

A binary string $f$ is a factor of string $u$ if $f$ appears as a sequence of $|f|$ consecutive bits of $u$, where $|f|$ denotes the length of $f$.
Ili\'{c} et al. \cite{asy1} introduced
 \emph{generalized Fibonacci cube}, $Q_{d}(f)$, as the graph obtained from $Q_{d}$ by
removing all strings that contain a given binary string $f$ as a factor.
In this notation the classical Fibonacci cube $\Gamma_{d}$ is the graph $Q_{d}(11)$.
The subclass $Q_{d}(1^{s})$ of generalized Fibonacci cube has been studied in \cite{whs,liu,nzs}.

A binary string $f$ is called good if $Q_{d}(f)$ is
an isometric subgraph of  $Q_{d}$ for all $d\geq1$, bad otherwise.
It was estimated that about eight percent of all strings are good \cite{sks}.
The index of a binary string $f$, denoted $B(f)$,
is the smallest integer $d$ such that $Q_{d}(f)$ is not an isometric subgraph of $Q_{d}$.
Ili\'{c} et al. \cite{asy} showed $B(f)<2|f|$ for almost all bad string $f$ and
thus posed the following conjecture:

\vspace{0.2cm}
\trou\noi{\bf Conjecture 1.1.(\cite{asy})}
\emph{For any bad string $f$, $B(f)<2|f|$.}

They also posed the following conjecture.

\vspace{0.2cm}
\trou\noi{\bf Conjecture 1.2.(\cite{asy1} )}
\emph{If $Q_{d}(f)$ is an isometric subgraph of $Q_{d}$,
then $Q_{d}(ff)$ is an isometric subgraph of $Q_{d}$.}

In this paper, our main aim is to confirm such conjectures.
The organization of the remainder of this paper is as follows.

In Section 2, we introduce some preliminary concepts and results of generalized Fibonacci cubes.
In Section 3, we obtain a basic result that if there exist $p$-critical words for $Q_{B(f)}(f)$,
then $p$=2 or $p$=3 by revealing the structures of $p$-critical words for $Q_{B(f)}(f)$.
In section 4, we prove that both Conjectures 1.1 and 1.2 to be true.

\section{Preliminaries}

In this section we introduce some concepts and results of generalized Fibonacci cube need in this paper.

For a binary string $f=f_{1}f_{2}\cdots f_{d}$ let $f^{R}=f_{d}\cdots f_{2}f_{1}$
be the reverse of $f$ and $\overline{f}=\overline{f}_{1} \overline{f}_{2}\cdots \overline{f}_{d}$ be the binary complement of $f$, where $\overline{f}_{i}=1-f_{i}$, $i=1,\ldots,d$.
Let $|f|$ denote the length of $f$.
Denote the string with 1 in coordinate $i$ and 0 elsewhere with $e_{i}$.
For strings $\alpha$ and $\gamma$ of equal length,
let $\alpha+\gamma$ denote their sum computed bitwise modulo 2.
In particular,
$\alpha+e_{i}$ is the string obtained from $\alpha$ by complementing its $i^{th}$ bit.


For a connected graph $G$,
the distance $d_{G}(\mu, \nu)$ between vertices $\mu$ and $\nu$ is the length of a shortest $\mu,\nu$-path.
Given two binary strings $\alpha$ and $\beta$ with the same length,
their Hamming distance $H(\alpha,\beta)$ is the number of bits in which they differ.
It is known that \cite{ik} for any vertices $\alpha$ and $\beta$ of $Q_{d}$, $H(\alpha,\beta)=d_{Q_{d}}(\alpha,\beta)$.

Obviously, for any subgraph $H$ of $G$, $d_{H}(\mu,\nu)\geq d_{G}(\mu,\nu)$.
If $d_{H}(\mu,\nu)=d_{G}(\mu,\nu)$ for all $\mu,\nu\in V(H)$,
then $H$ is called an \emph{isometric subgraph} of $G$,
and simply write $H\hookrightarrow G$, and $H\not\hookrightarrow G$ otherwise.
More generally, let $H$ and $G$ be any arbitrary graphs.
A mapping $h: V(H)\rightarrow V(G)$
is an \emph{isometric embedding} of $H$ into $G$ if
$d_{H}(\mu,\nu)=d_{G}(h(\mu),h(\nu))$ for any $\mu$ and $\nu\in V(H)$.

For two vertices $\mu$ and $\nu$ of graph $G$,
the set of vertices lying on shortest $\mu,\nu$-paths is called the interval between $\mu $
and $\nu$, denoted by $I_{G}(\mu,\nu)$.
Let $\alpha$ and $\beta\in Q_{d}(f)$ and $p\geq2$.
Then $\alpha$ and $\beta$ are called \emph{$p$-critical words} \cite{asy1} for $Q_{d}(f)$ if $d_{Q_{d}}(\alpha,\beta)=H(\alpha,\beta)=p$,
but none of the neighbors of $\alpha$ in $I_{Q_{d}}(\alpha,\beta)$
 belongs to $Q_{d}(f)$ or none of the neighbors of $\beta$ in
  $I_{Q_{d}}(\alpha,\beta)$ belongs to $Q_{d}(f)$.

The sufficiency of the following lemma was given by Lemma 2.4 of \cite{asy1}.

\trou\noi{\bf Lemma 2.1.} \emph{Suppose that $f$ is a nonempty binary string and $d\geq1$.
Then $Q_{d}(f)\not\hookrightarrow Q_{d}$ if and only if there exist
$p$-critical words for $Q_{d}(f)$ for some $p\geq2$.}

\trou\noi{\bf Proof}. We only show the necessity.
Assume that $Q_{d}(f)\not\hookrightarrow Q_{d}$,
then there exist vertices $\alpha$ and $\beta \in V(Q_{d}(f))$ such that $d_{Q_{d}(f)}(\alpha,\beta)>H(\alpha,\beta)=p\geq2$.
Let $\alpha$ and $\beta$ be such vertices with the minimum $p$.
Set $\alpha=a_{1}a_{2}\cdots a_{d}$, $\beta=b_{1}b_{2}\cdots b_{d}$,
$a_{i_{j}}\neq b_{i_{j}}$, and $\alpha_{j}$=$\alpha+e_{i_{j}}$, where $j=1,2,\ldots,p$.
We claim that $\alpha_{k}\not\in V(Q_{d}(f))$ for all $k\in\{1,2,\ldots,p\}$.
On the contrary suppose that $\alpha_{k}\in V(Q_{d}(f))$ for some $k\in\{1,2,\ldots,p\}$,
then $H(\alpha_{k},\beta)=p-1$.
If $p-1=1$ or $p-1=d_{Q_{d}(f)}(\alpha_{k},\beta)$,
then $p=H(\alpha,\beta)=d_{Q_{d}(f)}(\alpha,\beta)$,
a contradiction.
So $2\leq p-1=H(\alpha_{k},\beta)<d_{Q_{d}(f)}(\alpha_{k},\beta)$.
It is a contradiction to the minimality of $p$.
Hence the claim holds.
Thus $\alpha$ and $\beta$ are $p$-critical words for $Q_{d}(f)$. $\hfill\blacksquare$


The following lemma holds for a bad string.

\trou\noi{\bf Lemma 2.2 (\cite{sks})} \emph{Suppose that $Q_{d}(f )\not\hookrightarrow Q_{d}$. Then $Q_{d'}(f)\not\hookrightarrow Q_{d'}$ for all $d'>d$.}

Letting $f$ be any bad string,
$Q_{d}(f)\hookrightarrow Q_{d}$ for only a finite number of dimensions $d< B(f)$ by Lemma 2.2.

 Lemma 2.3 holds obviously.
It will be used repeatedly in proving the main results of next sections.

\trou\noi{\bf Lemma 2.3}. \emph{Let $f=f_{1}f_{2}\cdots f_{|f|}$ be a binary string,
$r,s $ and $t$ be positive integers such that $|f|\geq t+r+s$.
If any two equations of the three equations $(1)$ hold,
then the third must hold,}
\begin{equation}
\begin{split}
 f_{t}=f_{t+r}, f_{t+r}=f_{t+r+s},f_{t}=f_{t+r+s}.\\
\end{split}
\end{equation}


\trou\noi{\bf Lemma 2.4}.
\emph{Let $f=f_{1}f_{2}\cdots f_{|f|}$ be a binary
string and $r,s,t$ be positive integers such that $|f|> r+s-g$
and $t\leq |f|+g-r-s$, where $g$ is the greatest common divisor of $r$ and $s$.
If $\frac{s+r}{g}-1$ ones of the $\frac{s+r}{g}$ equations $(2)$ and $(3)$ hold,
then the remaining one must hold,}
\begin{equation}
\begin{split}
f_{t+(i-1)g}=f_{t+(i-1)g+s}, i=1,\ldots,\frac{r}{g},\\
\end{split}
\end{equation}
\begin{equation}
\begin{split}
f_{t+(j-1)g}=f_{t+(j-1)g+r}, j=1,\ldots,\frac{s}{g}.
\end{split}
\end{equation}

\trou\noi{\bf Proof}. Letting $k_{1}=\frac{r}{g}$ and $k_{2}=\frac{s}{g}$,
$k_{1}$ and $k_{2}$ are relatively prime.
We first construct the following 2-regular bipartite
 graph $G^{\langle r,s\rangle}$ with the bipartition $(X,Y)$,
where

$X=\{v^{1}_{i}~|~i=1,\ldots,k_{1}\} \cup\{v^{3}_{i}~|~i=1,\ldots,k_{2}\}$ and

$Y=\{v^{2}_{i}~|~i=1,2,\ldots,k_{1}+k_{2}\}$.

The edge set of $G^{\langle r,s\rangle}$ is

$E(G(t^{\langle r,s\rangle}))=\{v^{1}_{i}v^{2}_{i},v^{1}_{i}v^{2}_{i+k_{2}}~|~i=1,\ldots,k_{1}\}
\cup\{v^{3}_{k}v^{2}_{k},v^{3}_{k}v^{2}_{k+k_{1}}~|~k=1,\ldots,k_{2}\}$.

An example of $G^{\langle r,s\rangle}$ that $k_{1}=10$ and $k_{2}=3$ is shown in Fig . 1.

\begin{figure}[h]
\begin{center}
\includegraphics[scale=0.85]{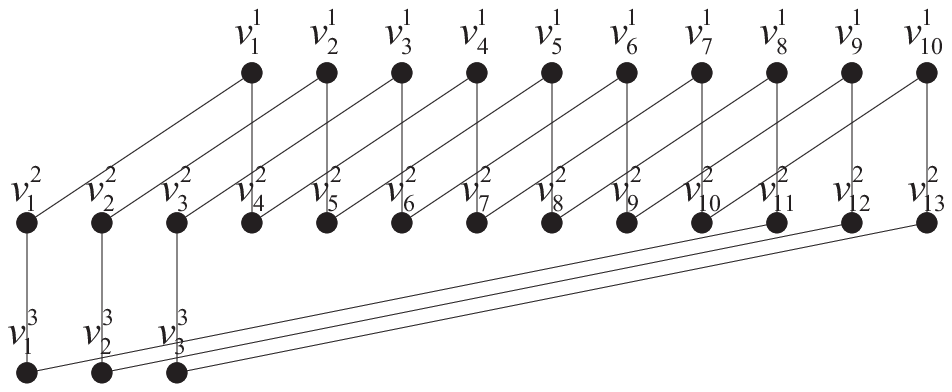}\\
{{\small Fig. 1. An example of the graph $G ^{\langle r,s\rangle}$ such that $k_{1}=10$ and $k_{2}=3$}}
\end{center}
\end{figure}

Letting $v^{1}_{i}$ denote $f_{t+(i-1)g}$ for $i=1,\ldots,k_{1}$,
$v^{2}_{j}$ denote $f_{t+(j-1)g}$ for $j=1,2,\ldots,k_{1}+k_{2}$
and $v^{3}_{k}$ denote $f_{t+(k-1)g+r}$ for $k=1,\ldots,k_{2}$,
a vertex in $X$ is joined to a vertex in $Y$ if they denote some one equation of (2), (3) and (4), where (4) are the following equations,

\begin{equation}
\begin{split}
f_{t+(i-1)g}=f_{t+(i-1)g}, i=1,\ldots,\frac{r+s}{g}.
\end{split}
\end{equation}

Now we show that $G^{\langle r,s\rangle}$ is connected.
Obviously if this holds,
then $G^{\langle r,s\rangle}$ is a cycle, and so this lemma is proved.
Since $v^{1}_{i}$ is joined to $v^{2}_{i+k_{2}}$ for $ i=1,\ldots,k_{1}$,
and $v^{3}_{k}$ is joined to $v^{2}_{k}$ for $k=1,\ldots,k_{2}$,
we only need to show that any two vertices of $Y$ are connected,
and this can be done by proving that
$v^{2}_{1}$ and $v^{2}_{j}$ are connected for $ j=1,\ldots,k_{1}+k_{2}$.
Two claims are needed.

\trou\noi{\bf Claim 1}.
Let $Z=\{h_{j}\  | \ h_{j}:=(j-1)k_{1}\mod (k_{1}+k_{2})$,
$0\leq h_{j} \leq k_{1}+k_{2}$ and $j=1,2,\ldots,k_{1}+k_{2}\}$.
Then $Z=\{0,1,\ldots,k_{1}+k_{2}-1\}$,
and $h_{j}=k_{2}$ if and only if $j=k_{1}+k_{2}$.

\noindent
{\em Proof}. First we show that $Z=\{0,1,\ldots,k_{1}+k_{2}-1\}$.
We only need to prove that if $1\leq q_{1}<q_{2}\leq k_{1}+k_{2}$,
then $h_{q_{2}}\neq h_{q_{1}}$.
Otherwise there exist $m\geq1$ such that $(q_{2}-q_{1})k_{1}=m(k_{1}+k_{2})$,
that is, $(q_{2}-q_{1}-m)k_{1}=mk_{2}$.
If $k_{1}=k_{2}$, then $k_{1}=k_{2}=1$ and $q_{2}-q_{1}=2m$,
but this is a contradiction to $1\leq q_{1}<q_{2}\leq 2$;
if $k_{1}\neq k_{2}$, then $q_{2}-q_{1}=k_{2}+m$ and $k_{1}=m$,
but this is a contradiction to $1\leq q_{1}<q_{2}\leq k_{1}+k_{2}$.
Hence $Z=\{0,1,\ldots,k_{1}+k_{2}-1\}$.

Now we show that $h_{j}=k_{2}$ if and only if $j=k_{1}+k_{2}$.
Since $Z=\{0,1,\ldots,k_{1}+k_{2}-1\}$, we only need to show the sufficiency.
Suppose that $(k_{1}+k_{2}-1)k_{1}=m(k_{1}+k_{2})+h_{j}$, where $m\geq0$.
Then $(k_{1}+k_{2})(k_{1}-m)= k_{1}+h_{j}$.
Obviously, $k_{1}-m\geq1$.
If $k_{1}-m\geq 2$,
then $ k_{1}+h_{j}=(k_{1}+k_{2})(k_{1}-m)\geq 2(k_{1}+k_{2})$,
it is a contradiction to $h_{j}< k_{1}+k_{2}$.
Hence $k_{1}-m=1$ and so $h_{j}=k_{2}$.
$\hfill\Box$

Recall that $v^{2}_{1}$ is joined to $v^{3}_{1}$,
and $v^{3}_{1}$ is joined to $v^{2}_{1+k_{1}}$ in graph $G^{\langle r,s\rangle}$.
Hence we get a path of $G^{\langle r,s\rangle}$ that start from the first vertex $v^{2}_{1}$ of $Y$ and reach the
second vertex $v^{2}_{1+k_{1}}$ of $Y$,
if we continue to walk along this path as far as possible, we have the following claim:

\trou\noi{\bf Claim 2}. The $k^{th}$ vertex of $Y$ that walking along this path is $v^{2}_{1+h_{k}}$, $k=1,2,\ldots,k_{1}+k_{2}$.

\noindent
{\em Proof}. Recall that $k_{1}$ and $k_{2}$ are relatively prime.
If $k_{1}=k_{2}$,
then $G^{\langle r,s\rangle}$ is 4-cycle,
and so the claim holds.

Now we prove this claim by induction on $k$  for the case $k_{1}\neq k_{2}$.
First we show it holds $k=1,2$.
In fact, the first vertex is $v^{2}_{1}=v^{2}_{1+h_{1}}$
and the second one is $v^{2}_{1+k_{1}}=v^{2}_{1+h_{2}}$.

We suppose that it holds for $2\leq k <k_{1}+k_{2}$ and $h_{k}=(k-1)k_{1}-m(k_{1}+k_{2})$ for some $m\geq 0$.
As shown in Claim 1, $h_{k}\neq k_{2}$, and so $h_{k}<k_{2}$ or $h_{k}>k_{2}$.
If $h_{k}<k_{2}$,
then $v^{2}_{1+h_{k}}$ is joined to $v^{3}_{1+h_{k}}$ and
$v^{3}_{1+h_{k}}$ is joined to $v^{2}_{1+h_{k}+k_{1}}$.
Now we show that $h_{k}+k_{1}=h_{1+k}$.
In fact, $kk_{1}=h_{k}+k_{1}+m(k_{1}+k_{2})$.
So $h_{k+1}= kk_{1} \mod (k_{1}+k_{2})=h_{k}+k_{1}$.
If $h_{k}>k_{2}$, then $v^{2}_{1+h_{k}}$ is joined to $v^{1}_{1+h_{k}-k_{2}}$ and
$v^{1}_{1+h_{k}-k_{2}}$ is joined to $v^{2}_{1+h_{k}-k_{2}}$.
Now we show that $h_{k}-k_{2}=h_{1+k}$.
In fact, $kk_{1}=h_{k}-k_{2}+(m+1)(k_{1}+k_{2})$ and so
$h_{k+1}\equiv kk_{1}($mod$(k_{1}+k_{2}))=h_{k}-k_{2}$.
So the claim holds for $k+1$.
Hence by inductive assumption the claim holds.
$\hfill\Box$

So by Claims 1 and 2,
we show that $v^{2}_{1}$ and $v^{2}_{j}$ are connected for $ j=1,\ldots,k_{1}+k_{2}$.
This completes the proof.
$\hfill\blacksquare$


Now we consider a special case of Lemma 2.4 such that $r+s=|f|$.

\trou\noi{\bf Corollary 2.5}.
\emph{Let $f=f_{1}f_{2}\cdots f_{|f|}$ be a binary string and $r,s$
be positive integers such that $|f|=r+s$.
Suppose $f_{i}=f_{i+s}$ for $i=1,\ldots,r$ and $f_{j}=f_{j+r}$ for $j=1,\ldots,s$,
then equations $(2)$ and $(3)$ holds for $t=1,\ldots,g$.}

By the structure of the graph $G^{\langle r,s \rangle}$ of Lemma 2.4,
we have the following lemma.

\trou\noi{\bf Lemma 2.6}. \emph{Let $r,s,t,$ and $g$ be the same as in Lemma $2.4.$}

$(\romannumeral1)$ If equations (2) hold for all $i\in \{1,\ldots,\frac{r}{g}\}\setminus \{i_{1},i_{2}\}$ and equations (3) hold for all $ j\in \{1,\ldots,\frac{s}{g}\}$,
then $f_{t+(i_{1}-1)g}=f_{t+(i_{2}-1)g+s}$ and $f_{t+(i_{2}-1)g}=f_{t+(i_{1}-1)g+s}$.

$(\romannumeral2)$ If equations (2) hold for all $i\in \{1,\ldots,\frac{r}{g}\}$ and equations (3) hold for all $ j\in \{1,\ldots,\frac{s}{g}\}\setminus \{j_{1},j_{2}\}$,
then $f_{t+(j_{1}-1)g}=f_{t+(j_{2}-1)g+r}$ and $f_{t+(j_{2}-1)g}=f_{t+(j_{1}-1)g+r}$.

$(\romannumeral3)$ If equations (2) hold for all $i\in \{1,\ldots,\frac{r}{g}\}\setminus \{i_{1}\}$ and equations (3) hold for all $ j\in \{1,\ldots,\frac{s}{g}\}\setminus \{j_{1}\}$, then $f_{t+(i_{1}-1)g}=f_{t+(j_{1}-1)g}$ and $f_{t+(i_{1}-1)g+s}=f_{t+(j_{1}-1)g+r}$.

\section{A basic result on $p$-critical words}

In this section we prove that for every bad string $f$,
the $p$-critical words of $Q_{B(f)}(f)$ satisfy that $p=2$ or $p=3$.
More generally, there exist $2$-critical words or $3$-critical words for $Q_{d}(f)$,
where $d \geq B(f)$.

Suppose that $f$ is a binary string, $\alpha\in V(Q_{d}(f))$ and $\alpha'=\alpha+e_{j} \not\in V(Q_{d}(f))$ for some $j\in\{1,2,\ldots,d\}$,
in other words, $\alpha'$ contains the factor $f$.
We denote $f^{(u_{j})}$ as the factor $f$ appears in $\alpha'$ starting from coordinate $u_{j}$. 
Note that if $\alpha''=\alpha+e_{k} \not\in V(Q_{d}(f))$ for some $k\in\{1,2,\ldots,d\}\setminus\{j\}$, then $u_{j}\neq u_{k}$.

In the lemmas and theorems of this section, we suppose that
$f=f_{1}f_{2}\cdots f_{|f|}$ is a bad string and $d_{0}=B(f)$;
$\alpha=a_{1}a_{2}\cdots a_{d_{0}}$ and
$\beta=b_{1}b_{2}\cdots b_{d_{0}}$ are $p$-critical
words for $Q_{d_{0}}(f)$, $a_{i_{j}}\neq b_{i_{j}}$ and $\alpha_{j}=\alpha+e_{i_{j}}\not\in V(Q_{d_{0}}(f))$ for $j=1,2,\ldots,p$.

\trou\noi{\bf Lemma 3.1}. \emph{For any $k\in\{1,2,\ldots,p\}$,
there exists some one $j \in\{1,2,\ldots,p\} \setminus \{k\}$ such that
$i_{k}\in [u_{i_{j}},u_{i_{j}}+|f|-1]$}.

\trou\noi{\bf Proof}.
We distinguish three cases: $p=2$, $p=3$ and $p\geq 4$.

Obviously if $p=2$, then $i_{1}\in [u_{i_{2}},u_{i_{2}}+|f|-1]$ and $i_{2}\in [u_{i_{1}},u_{i_{1}}+|f|-1]$, otherwise $\beta$ contains the factor $f$.

Let $p=3$. On the contrary we distinguish the following three cases:

$i_{1}\not\in [u_{i_{2}},u_{i_{2}}+|f|-1]$
and $i_{1}\not\in [u_{i_{3}},u_{i_{3}}+|f|-1]$,

$i_{2}\not\in [u_{i_{1}},u_{i_{1}}+|f|-1]$
and~$i_{2}\not\in [u_{i_{3}},u_{i_{3}}+|f|-1]$ and

$i_{3}\not\in [u_{i_{1}},u_{i_{1}}+|f|-1]$
and $i_{3}\not\in [u_{i_{2}},u_{i_{2}}+|f|-1]$.

We only show that the first and second cases do not hold since the third one can
be proved similarly as the first one.
In the first case, we claim that $i_{2}\in [u_{i_{3}},u_{i_{3}}+|f|-1]$
and $i_{3}\in [u_{i_{2}},u_{i_{2}}+|f|-1]$,
otherwise $\beta $ contains $f$ as a factor, a contradiction.
Let $\alpha'$ and $\beta'$ be the strings obtained by deleting
the first to $i_{1}^{th}$ bit of $\alpha$ and $\beta$, respectively.
Then $\alpha'$ and $\beta'$ are $2$-critical words for $Q_{d'}(f)$, where $d'=d_{0}-i_{1}$,
it is a contradiction to $d_{0}=B(f)$.
In the second case, obviously the factor that starts from the
$ u_{i_{1}}^{th}$ to $(u_{i_{1}}+|f|-1)^{th}$ bit of $\beta$ is $f$, it is a contradiction.

Now we consider $p\geq 4$.
Assume for some one $k\in\{1,2,\ldots,p\}$,
$i_{k}\not\in [u_{i_{j}},u_{i_{j}}+|f|-1]$ holds for
all $j\in \{1,2,\ldots,p\}\setminus \{ k\}$.
Then $u_{i_{j}}+|f|-1<i_{k}$ for all $j\in X=\{1,\ldots,k-1\}$
and $u_{i_{j}}>i_{k}$ for all $j\in Y=\{k+1,\ldots,p\}$.
Note that $X=\emptyset$ if $k=1$, and $Y=\emptyset$ if $k=p$.
Since $p\geq 4$, at least one of $|X|$ and $|Y|$ not less than 2.
Without loss of generality, suppose $|Y|\geq2$.
Letting $\alpha'$ and $\beta'$ be the strings obtained by deleting the
first bit to the $i_{k}^{th}$ bit of $\alpha$ and $\beta$ respectively,
then $\alpha'$ and $\beta'$ are $p'$-critical words for $Q_{d'}(f)$,
where $p'=p-k$ and $d'=d_{0}-i_{k}$, it is a contradiction to that $d_{0}=B(f)$.
$\hfill\blacksquare$

\trou\noi{\bf Theorem 3.2}. \emph{If there exist $p$-critical words for $Q_{d_{0}}(f)$,
then $p=2$ or $p=3$.}

\trou\noi{\bf Proof}. Letting $k=1$, there must exist $j\in\{2,\ldots,p\}$
such that $i_{1}\in [u_{i_{j}},u_{i_{j}}+|f|-1]$ by Lemma 3.1.
Suppose that $j$ be the smallest one among all such coordinates.
According to $i_{j}$ belongs to $[u_{i_{1}}, u_{i_{1}}+|f|-1]$ or not,
we distinguish two cases.

\textbf{Case 1.} $i_{j}\in [u_{i_{1}},u_{i_{1}}+|f|-1]$.

There are two subcases: $u_{i_{1}}<u_{i_{j}}$ or $u_{i_{1}}>u_{i_{j}}$.

\textbf{Subcase 1.1.} $u_{i_{1}}<u_{i_{j}}$.

Let $r=u_{i_{j}}-u_{i_{1}}$. Obviously, $|f|-2\geq r\geq1$.
Let $l=i_{1}-u_{i_{1}}+1$ and $m=i_{j}-u_{i_{1}}+1$.
Set $\alpha'=a'_{1}\cdots a'_{r+|f|}=a _{u_{i_{1}}}\cdots a _{u_{i_{j}}+|f|-1}$.
Then $a'_{l}=a_{i_{1}}$ and $a'_{m}=a_{i_{j}}$.
It follows that $\alpha'+e_{l}$ and $\alpha'+e_{m}$ both contain the factor $f$ .
Obviously the first coordinate of $f^{(u_{l})}$ and $f^{(u_{m})}$
are $u_{l}=1$ and $u_{m}=r+1$, respectively.
Since $\alpha'$ is a factor of $\alpha$, $\alpha'$ does not contain the factor $f$.
Let $\beta'=\alpha'+e_{l}+e_{m}$.
We claim that $\beta'$ contains no factor $f$.
Obviously it holds for $r=1$.
On the contrary suppose that $\beta'$ contains factor $f$ for $r\geq2$,
and denote this factor as $f'$.
Then the first coordinate $v$ of  $f'$ must satisfy with $1\leq v \leq r$ by $|\beta'|=r+|f|\leq 2|f|-2$.
So this factor must contain  $\overline{a'}_{l}$ and $\overline{a'}_{m}$ of $\beta'$.
Set $s'=r+1-v$ and $r'=v-1$.
Hence there are three copies of factor $f$:
$f^{(r+1)}$, $f'$ and $f^{(1)}$.
They have the relative positions as shown in Fig. 2 and satisfy that:

$f^{(r+1)}_{t}=f'_{t}=f^{(1)}_{t}$ for all $t\in\{1,\ldots,|f|\}$;

$f^{(r+1)}_{t}=f'_{t+s'}=a'_{t+r}$ for all $t\in\{1,\ldots,|f|-s'\}\setminus \{l-r\}$,

specially, $f^{(r+1)}_{m-r}=f'_{m-r'}=\overline{a'}_{m}$, but
$f^{(r+1)}_{l-r}=a'_{l}$, $f'_{l-r'}=\overline{a'}_{l}$;

$f'_{t}=f^{(1)}_{t+r'}$, $t\in\{1,\ldots,|f|-r'\} \setminus \{m-r'\}$,

specially, $f'_{l-r'}=f^{(1)}_{l}=\overline{a'}_{l}$,
but $f'_{m-r'}=\overline{a'}_{m}$, $f^{(1)}_{m}=a'_{m}$.

\begin{figure}[h]
\begin{center}
\includegraphics[scale=0.85]{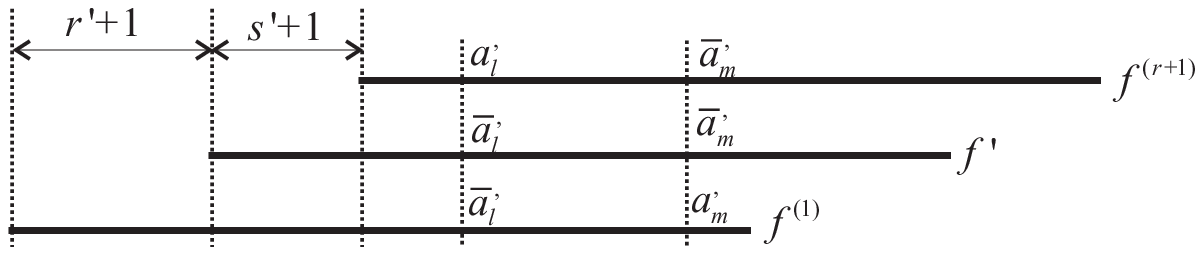}\\
{{\small Fig. 2. Illustration of the claim that $\beta'$
 contains no factor $f$ in Subcase 1.1.}}
\end{center}
\end{figure}

Obviously $f'_{l-r}=f^{1}_{l-s'}, f^{(r+1)}_{l-s'}=f'_{l}$,
$a'_{l}=f^{(r+1)}_{l-r}$ and $ f^{(1)}_{l}=\overline{a'}_{l}$,
in other words, equations (1) holds for $t=l-r$.
By Lemma 2.3, $a'_{l}=\overline{a'}_{l}$. It is impossible.
Thus $\beta'$ contains no $f$.

Hence $\alpha'$ and $\beta'$ are 2-critical words for $Q_{d'}(f)$,
where $d'=r+|f|$ and $1\leq r \leq |f|-2$.
So $\alpha=\alpha'$, $\beta=\beta'$ and $B(f)=d'=r+|f|$ by the fact that $\alpha'$ is a factor of $\alpha$ and $d_{0}=B(f)$.
Thus $p=2$ in this subcase.

\textbf{Subcase 1.2.} $u_{i_{j}}<u_{i_{1}}$.

Letting $r=u_{i_{1}}-u_{i_{j}}$, $|f|-2\geq r\geq1$.
Set $l=i_{1}+1-u_{i_{j}}$ and $m=i_{j}+1-u_{i_{j}}$.
Denote $\alpha'=a' _{1}\cdots a' _{r+|f|}=a _{u_{i_{j}}}\cdots a _{u_{i_{1}}+|f|-1}$.
So $a'_{l}=a_{i_{1}}$ and $a'_{m}=a_{i_{j}}$.
It follows that $\alpha'+e'_{l}$ and $\alpha'+e'_{m}$ both contain factor $f$.
Obviously, the first coordinate of $f^{(u_{l})}$ and $f^{(u_{m})}$ are
$u_{l}=r+1$ and $u_{m}=1$, respectively.
Since $\alpha'$ is a factor of $\alpha$, $\alpha'$ contains no factor $f$.
Let $\beta'=\alpha'+e_{l}+e_{m}$.
If we show that $f$ is not a factor of $\beta'$,
then $\alpha=\alpha'$, $\beta=\beta'$ and $d'=r+|f|=B(f)$.
So $\alpha$ and $\beta$ are 2-critical words for $Q_{d}(f)$.
Thus $p=2$.

Now we assume that $\beta'$ contains factor $f$ and denote it as $f'$, then $r\geq2$.
Suppose the first coordinate of $f'$ is $v$.
Then $1<v\leq r$ and $f'$ must contain the $l^{th}$ and $m^{th}$ bits.
Of course $\alpha' \neq \alpha$, $\beta\neq\beta'$ and $p\geq3$,
in other words, there exist $k\in\{1,2,\ldots,p\}\setminus\{1,j\}$ such that
$i_{k}\in (i_{1},u_{i_{j}}+v+|f|-2]$.
In order to determine the coordinate of $i_{k}$, we need the following claims.

\textbf{Claim 1.} \emph{$r$ is even, $m=l+\frac{r}{2}$  and $\overline{a'}_{l}=a'_{m}$.}

\noindent
{\em Proof}.
Letting $r_{1}=v-1$ and $r_{2}=r+1-v$,
there are three copies of $f$:
$f^{(r+1)}$, $f'$ and $f^{(1)}$.
They have the relative positions as shown in Fig. 3 and satisfy that:

$f^{(r+1)}_{t}=f'_{t}=f^{(1)}_{t}=f_{t}$, $t\in\{1,\ldots,|f|\}$;

$f^{(r+1)}_{t}=f'_{t+r_{2}}=a'_{t+r}$, $t\in\{1,\ldots,|f|-r_{2}\} \setminus \{m-r\}$,

specially, $f^{(r+1)}_{l-r}=f^{(1)}_{l-r_{1}}=\overline{a'}_{l}$, but
$f^{(r+1)}_{m-r}=a'_{m}$, $f'_{m-r_{1}}=\overline{a'}_{m}$;

$f'_{t}=f^{(1)}_{t+r_{1}}=a'_{t+r_{1}}$, $t\in\{1,\ldots,|f|-r_{1}\}\setminus\{l-r_{1}\}$;

specially, $f'_{m-r_{1}}=f^{(1)}_{m}=\overline{a'}_{m}$, but
$f'_{l-r_{1}}=\overline{a'}_{l}$, $f^{(1)}_{l}=a'_{l}$.

\begin{figure}[h]
\begin{center}
\includegraphics[scale=0.85]{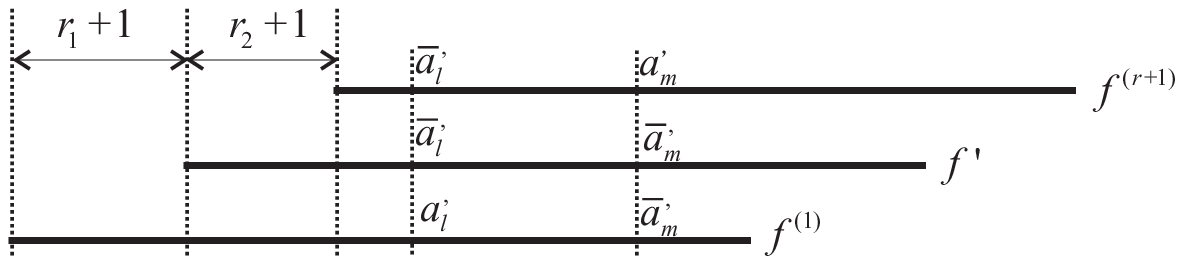}\\
{{\small Fig. 3. Illustration of the Claim 1 of Subcase 1.2.}}
\end{center}
\end{figure}

If $l+r_{1}\neq m$,
then equations (1) hold for $t=l-r$ and by Lemma 2.3,  $\overline{a'}_{l}=f^{(r+1)}_{l-r}=f^{(1)}_{l}=a'_{l}$.
Obviously it is impossible.
Hence $l+r_{1}=m$.
If $m-r_{2}\neq l$,
then equations (1) hold for $t=m-r$ and by Lemma 2.3,
$a'_{m}=f^{(r+1)}_{m-r}=f^{(1)}_{m}=\overline{a'}_{m}$.
Obviously it is impossible.
Hence $l+r_{2}=m$.

Thus $r_{1}=r_{2}=\frac{r}{2}$ and so $r=2r'$ is even.
Hence $\overline{a'}_{l}=f^{(r+1)}_{l-r}=f'_{l-r}=f^{(1)}_{l-r'}=
f^{(r+1)}_{l-r'}=f^{(r+1)}_{m-r}=a'_{m}$.
$\hfill\Box$

By Claim 1, $i_{j}-i_{1}=r/2=r'$ and  $\overline{a}_{i_{1}}=a_{i_{j}}$.
Note that the forbidden factor $f'$ contained in $\beta'$ means that $f$
is a factor of $\alpha+e_{i_{1}}+e_{i_{j}}$ such that starts from the bit $u_{i_{j}}+r'$,
and denote this factor as $f^{*}$.

By the above discussion, we know there must
exist $k\in\{1,2,\ldots,p\}\setminus \{1,j\}$ such that $ i_{k}\in ( i_{1},u_{i_{j}}+r'+|f|-1]$.
In the following claims we determine the position of $ i_{k}$ further.

\textbf{Claim 2.} \emph{There exists no $k$ such that $ i_{1}<i_{k}<i_{j}$.
In other words, $ j=2$.}

On the contrary suppose that there exists $ i_{k}$ such that $ i_{1}<i_{k}<i_{j}$.
Since $j$ is the minimal one such that $i_{1}\in [u_{i_{j}},u_{i_{j}}+|f|-1]$,
$u_{i_{k}}> i_{1}$.

\begin{figure}[h]
\begin{center}
\includegraphics[scale=0.85]{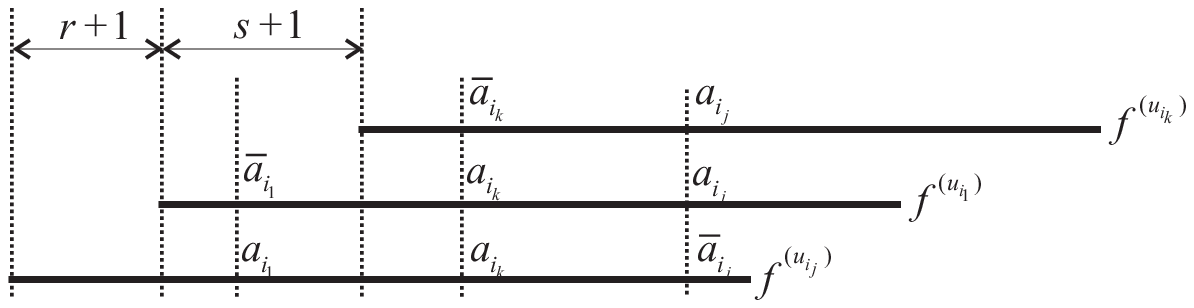}\\
{{\small Fig. 4. Illustration of the proof of Claim 2 of Subcase 1.2.}}
\end{center}
\end{figure}

Now we consider $f^{(u_{i_{k}})}$, $f^{(u_{i_{1}})}$ and $f^{(u_{i_{j}})}$.
Letting $s=u_{i_{k}}-u_{i_{1}}$,
those three copies have the relative positions as shown in Fig. 4 and they satisfy:

$f^{(u_{i_{k}})}_{t}=f^{(u_{i_{1}})}_{t}=f^{(u_{i_{j}})}_{t}$, $t\in\{1,\ldots,|f|\}$;

$f^{(u_{i_{k}})}_{t}=f^{(u_{i_{1}})}_{t+s}=a_{t+u_{i_{k}}-1}$, $t\in\{1, \ldots,|f|-s\} \setminus\{ i_{k}-u_{i_{k}}+1\}$,

$f^{(u_{i_{k}})}_{i_{k}-u_{i_{k}}+1}=\overline{a}_{i_{k}}$,
~$f^{(u_{i_{1}})}_{i_{k}-u_{i_{k}}+1+s}=a_{i_{k}}$;

$f^{(u_{i_{1}})}_{t}=f^{(u_{i_{j}})}_{t+r}=a_{t+u_{i_{1}}-1}$,
$t\in\{1, \ldots,|f|-r\} \setminus\{ i_{1}-u_{i_{1}}+1,~i_{j}-u_{i_{j}}+1\}$,

$f^{(u_{i_{1}})}_{ i_{1}-u_{i_{1}}+1}=\overline{a}_{i_{1}}$,
$f^{(u_{i_{j}})}_{ i_{1}-u_{i_{1}}+1+r}=a_{i_{1}}$,

$f^{(u_{i_{1}})}_{ i_{j}-u_{i_{j}}+1}=a_{i_{j}}$ and
$f^{(u_{i_{j}})}_{ i_{j}-u_{i_{j}}+1+r}=\overline{a}_{i_{j}}$.

If $i_{k}-i_{1}\neq s$,
then the equations (1) hold for $t=i_{k}-u_{i_{k}}+1$ and so by Lemma 2.3,
$a_{i_{k}}=\overline{a}_{i_{k}}$, a contradiction.
If $i_{k}-i_{1}=s$, then $i_{j}-i_{1}\neq s$,
and so the equations (1) hold for $t=i_{j}-u_{i_{k}}+1$ and so by Lemma 2.3,
$a_{i_{j}}=\overline{a}_{i_{j}}$, a contradiction.

Hence there exists no $k$ such that $ i_{1}<i_{k}<i_{j}$. In other words, $ j=2$.
$\hfill\Box$

By Claim 2, we get $i_{k}\in [i_{2}+1,u_{i_{2}}+r'+|f|-1]$.

\textbf{Claim 3.} \emph{There exists no $k$ such that $u_{i_{k}}\leq i_{1}$.}


We first prove that there exist no $k$ such that $u_{i_{k}}< u_{i_{2}}$.
Otherwise there exist $k$ such that $ i_{k}< u_{i_{2}}+|f|-1$.
Now we consider $f^{(u_{i_{1}})}$, $f^{(u_{i_{2}})}$ and $f^{(u_{i_{k}})}$.
They have the positions as shown in Fig. 5 $(a)$, where $s= u_{i_{2}}-u_{i_{k}}$,
and they satisfy:

$f^{(u_{i_{1}})}_{t}=f^{(u_{i_{2}})}_{t+r}
=a_{t+u_{i_{1}}-1}$, $i_{2}-u_{i_{1}}+1<t\leq |f|-r$,

$f^{(u_{i_{2}})}_{t}=f^{(u_{i_{k}})}_{t+s}
=a_{t+u_{i_{2}}-1}$, $1\leq t<i_{2}-u_{i_{2}}+1$.

\begin{figure}[h]
\begin{center}
\includegraphics[scale=0.85]{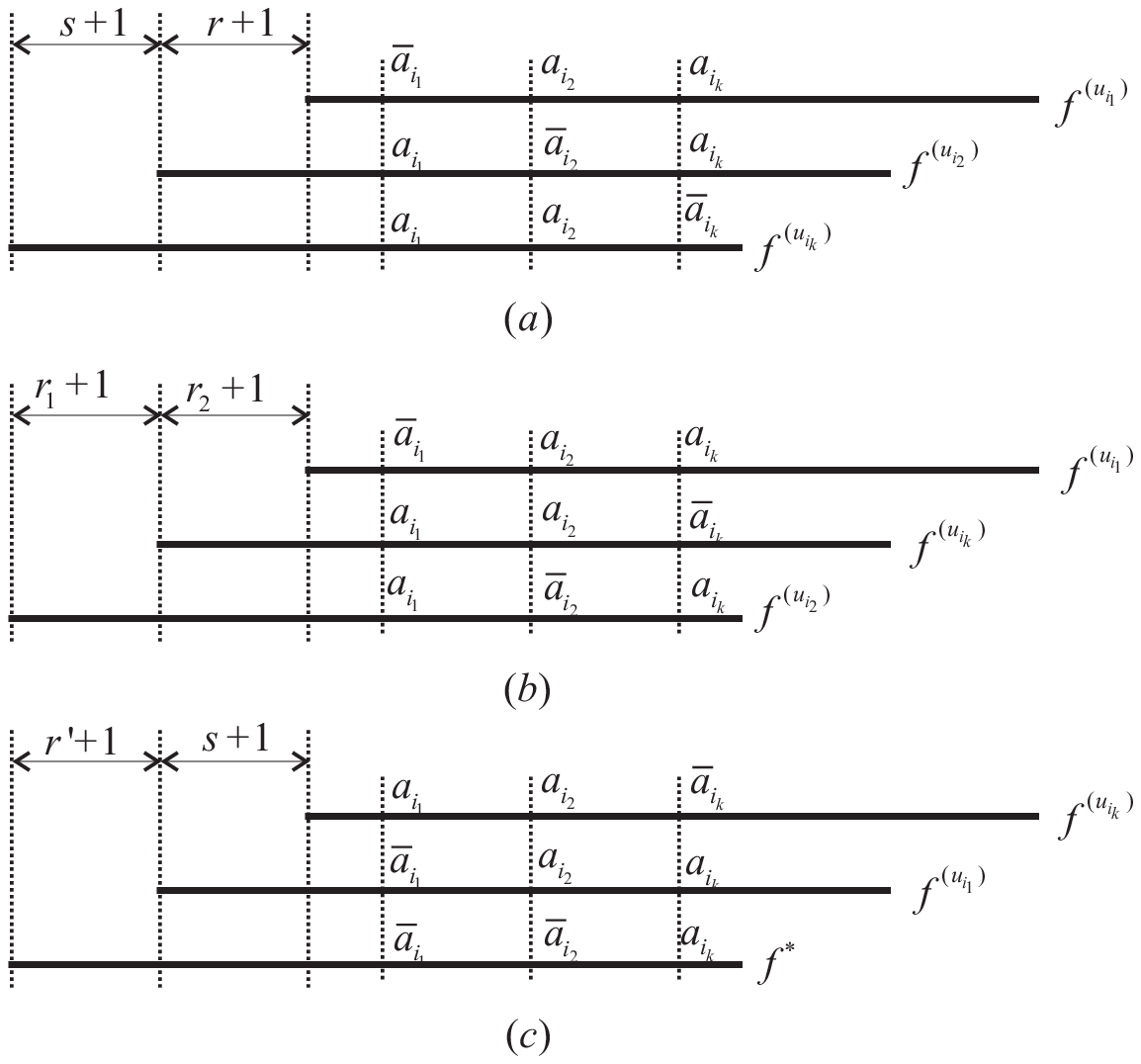}\\
{{\small Fig. 5. Illustration of the Claim 3 of Subcase 1.2.}}
\end{center}
\end{figure}

Note that $i_{2}-u_{i_{2}}>r$, $u_{i_{2}}+|f|-i_{2}-1>s$.
Hence the equations (2) and (3) hold for $t=i_{2}-u_{i_{1}}+1$,
and by Lemma 2.4,
$a_{i_{2}}=f_{u_{i_{2}}+|f|-i_{2}-1}=f_{u_{i_{2}}+|f|-i_{2}-1+r}=\overline{a}_{i_{2}}$.
Obviously it is impossible.
Thus there exists no $k$ such that $u_{i_{k}}< u_{i_{2}}$.

Next we prove that there exists no $ k$ such that $u_{i_{2}}< u_{i_{k}}< u_{i_{1}}$.
Otherwise there exists $ k$ such that $ i_{k}< u_{i_{1}}+|f|-1$.
Letting $u_{i_{k}}-u_{i_{2}}= r_{1}$ and $u_{i_{1}}-u_{i_{k}}=r_{2}$,
$r_{1}+r_{2}=r$ and $|f|\geq r_{1}+r_{2}+2$.
We consider $f^{(u_{i_{1}})}$,
$f^{(u_{i_{k}})}$ and $f^{(u_{i_{2}})}$ as shown in Fig. 5 $(b)$.
Note that $i_{2}-u_{i_{2}}> r_{1}+r_{2}+1=r+1$.
If $i_{k}-i_{1}\neq r_{1}$,
then equations (1) hold for $t=i_{1}-u_{i_{1}}+1$,
and so by Lemma 2.3, $\overline{a}_{i_{1}}=a_{i_{1}}$, a contradiction.
If $i_{k}-i_{1}=r_{1}$, then $i_{k}-i_{2}\neq r_{1}$.
Thus equations (1) hold for $t=i_{2}-u_{i_{1}}+1$
and so $a_{i_{2}}=\overline{a}_{i_{2}}$. It is impossible.
Thus there exists no $ k$ such that $u_{i_{2}}< u_{i_{k}}< u_{i_{1}}$.

Finally, we show that there exists no $ k$ such that $u_{i_{1}}< u_{i_{k}}\leq i_{1}$.
Otherwise we consider $f^{(u_{i_{k}})}$, $f^{(u_{i_{1}})}$ and $f^{*}$,
as it is shown in Fig. 5 $(c)$.
Recall that $f^{*}$ is the factor contained in $\alpha+e_{i_{1}}+e_{i_{j}}$ that starts from the $(v')^{th}$ bit, where $v'=u_{i_{2}}+r'$.
Note that $i_{2}-i_{1}=r'$ and equations (1) hold for $t=i_{1}-u_{i_{k}}+1$,
and so $a_{i_{1}}=\overline{a}_{i_{1}}$ by Lemma 2.3.
It is impossible.
This completes the proof of this claim.
$\hfill\Box$

\textbf{Claim 4.} \emph{There exists $k$ such that $i_{1}<u_{i_{k}}\leq i_{2}$,
and for such $k$, $i_{k}=i_{2}+r'$, $u_{ i_{k}}=u_{i_{1}}+r'$ and $a_{i_{k}}=\overline{a}_{i_{2}}$.}

\begin{figure}[h]
\begin{center}
\includegraphics[scale=0.80]{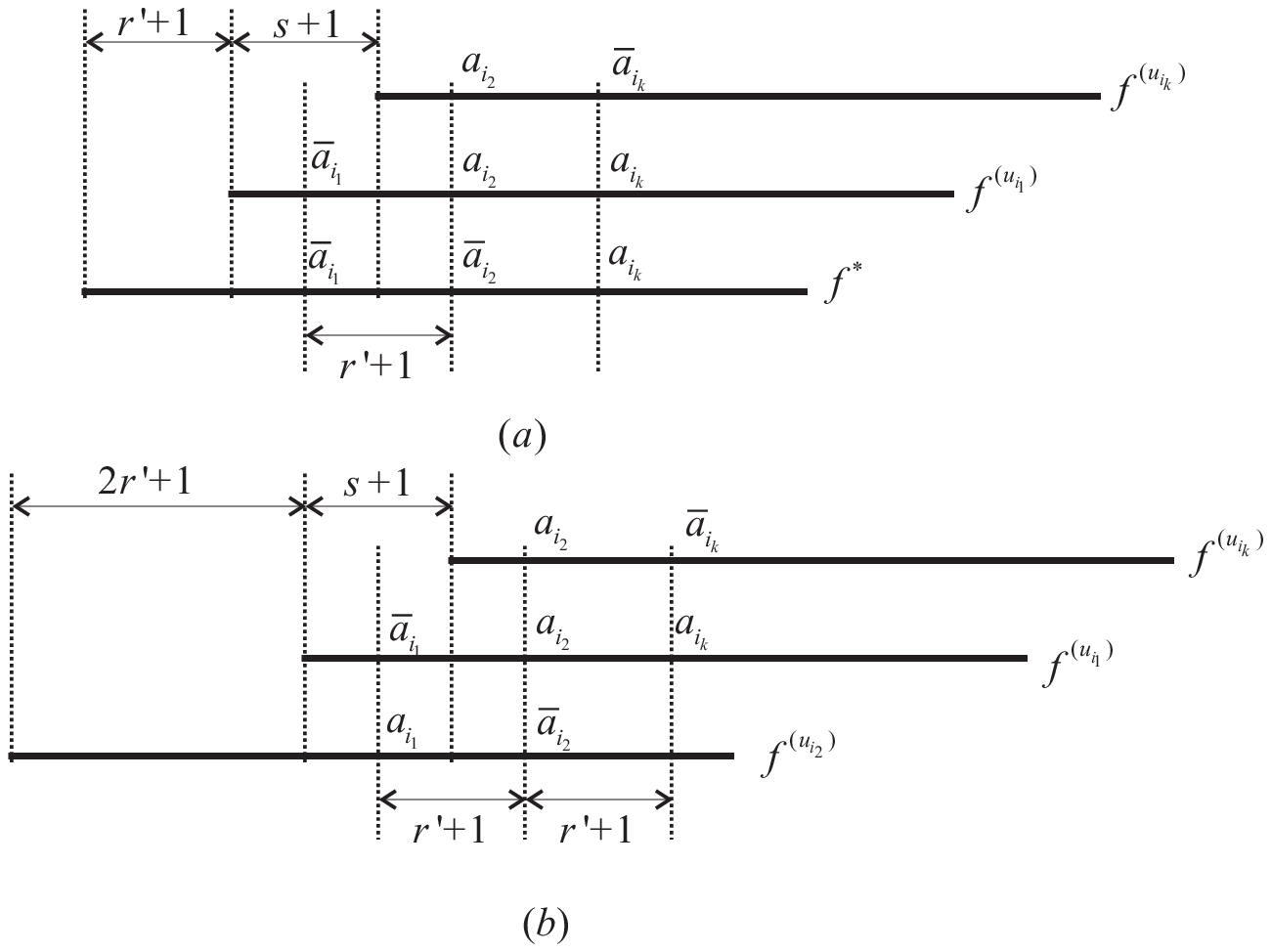}\\
{{\small Fig. 6. Illustration of the proof of Claim 4 of Subcase 1.2.}}
\end{center}
\end{figure}

By Claim 3,  $u_{i_{k}}> i_{1}$.
Now we show that there exists $u_{i_{k}}\leq i_{2}$.
Otherwise $ i_{2}<u_{i_{k}}$ for all $k\in \{3,\cdots,p\}$.
If $p=3$, then $k=3$.
So $\beta$ contains the factor $f^{(u_{i_{3}})}$, a contradiction.
Hence $p\geq4$. Letting $\mu$ and $\nu$ be the strings that obtained by deleting
the factor that start from the first bit to the $i_{2}^{th}$ bit of $\alpha$ and $\beta$, respectively,
$\mu$ and $\nu$ are $(p-2)$-critical words for $Q_{d'}(f)$,
where $d'=d_{0}-i_{2}$. It is a contradiction to $d_{0}=B(f)$.
Hence there exists $ k$ such that $u_{i_{k}}\leq i_{2}$.
Thus $i_{1}<u_{i_{k}}\leq i_{2}$.

Now we chose any $k$ such that $ i_{1}<u_{i_{k}}\leq i_{2}$.
First we show $ i_{k}=i_{2}+r'$.
Otherwise we consider $f^{(u_{i_{k}})}$, $f^{(u_{i_{1}})}$ and $f^{*}$,
as shown in Fig. 6 $(a)$. Hence equations (1) hold for $t=i_{2}-u_{i_{k}}+1$ and by Lemma 2.3, $a_{i_{2}}=\overline{a}_{i_{2}}$. obviously it is impossible.
Hence $ i_{k}=i_{2}+r'$ and so $a_{i_{k}}=\overline{a}_{i_{2}}$.
By $ i_{k}=i_{2}+r'$, $ i_{k}\in [i_{2}+1,u_{i_{2}}+r'+|f|-1]$.

Finally, we show that $u_{ i_{k}}=u_{i_{1}}+r'$.
Otherwise we consider $f^{(u_{i_{k}})}$, $f^{(u_{i_{1}})}$ and $f^{(u_{i_{2}})}$,
as shown in Fig. 6 $(b)$.
Hence equations (1) holds for $t=i_{2}-u_{ i_{k}}+1$.
By $i_{2}-i_{1}=i_{k}-i_{2}=r'$ and Lemma 2.3, $a_{i_{2}}=\overline{a}_{i_{2}}$.
It is impossible.
Thus  $u_{ i_{k}}=u_{i_{1}}+r'$.
$\hfill\Box$

By Claims 1---4, we find that $k\in\{1,2,\ldots,p\}\setminus\{1,j\}$ such that
$i_{k}\in (i_{1},v+|f|+1]$ is unique and $i_{k}=i_{1}+2r'=i_{2}+r'$,
where $3r'+1\leq |f|$.

Now we construct string $\alpha''=a'_{1}a'_{2}\cdots a'_{d'}
=a_{u_{i_{2}}}a_{u_{i_{2}}+1}\cdots a_{u_{i_{k}}+|f|-1}$.
Let $l =i_{1}-u_{i_{2}}+1$, $m=i_{2}-u_{i_{2}}+1$ and $n =i_{k}-u_{i_{2}}+1$,
then $a'_{l}=a_{i_{1}}, a'_{m}=a_{i_{2}}$ and $ a'_{n}=a_{i_{k}}$.
By Claims 1---4, $a'_{l}=\overline{a'}_{m}=a'_{n}$ and $n-m=m-l=r'$.
It follows that $\alpha''+e'_{l}$,
$\alpha''+e'_{m}$ and $\alpha''+e'_{l}$ all contain $f$ as factor,
and $u_{l}=2r'+1$, $u_{m}=1$ and $u_{n}=3r'+1$.
Obviously $\alpha''$ contains no $f$ as a factor since it is a factor of $\alpha$.

\begin{figure}[h]
\begin{center}
\includegraphics[scale=0.85]{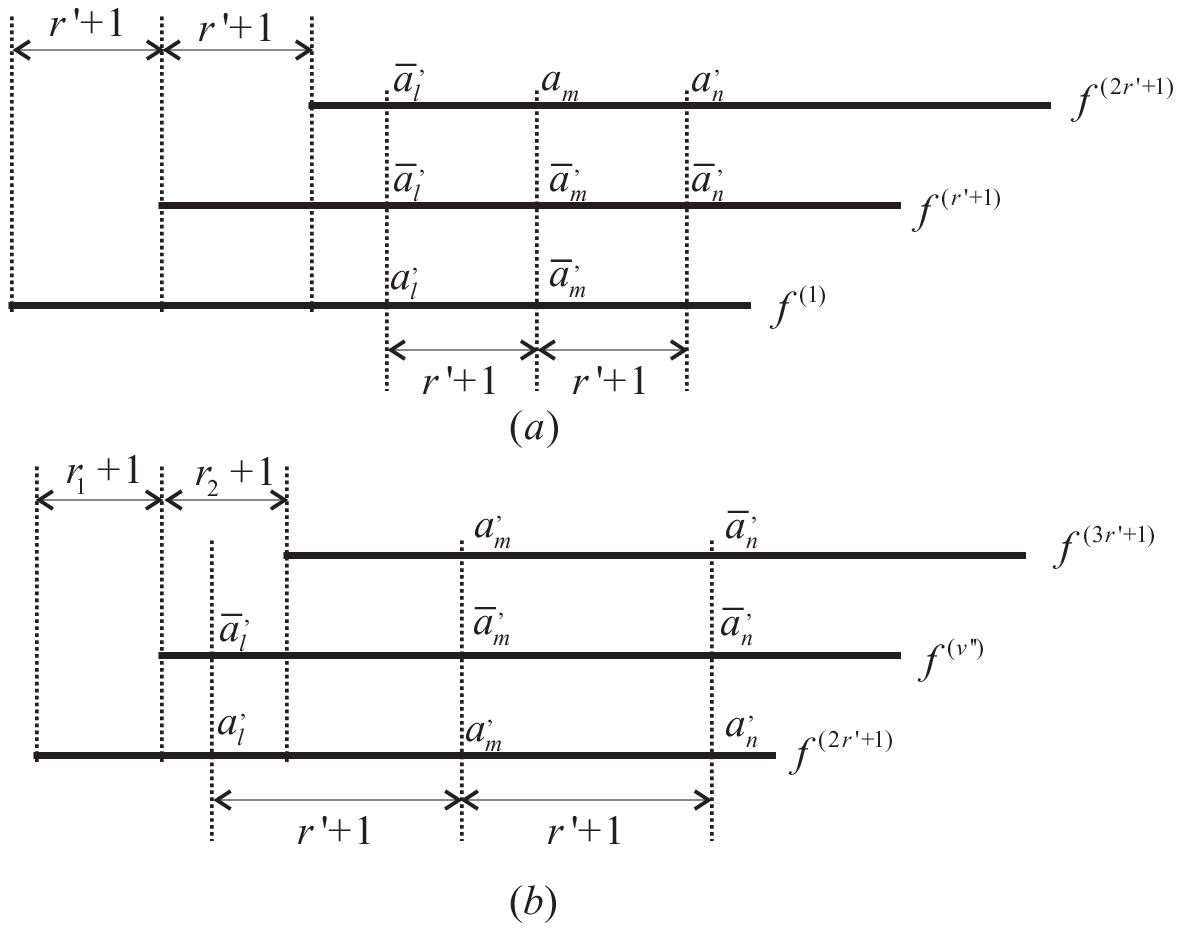}\\
{{\small Fig. 7. Illustration of the proof of $\beta''$ contains no $f$ as factor.}}
\end{center}
\end{figure}

Let $\beta''=\alpha''+e_{l}+e_{m}+e_{n}$.
Now we show that $\beta''$ does not contain $f$ as a factor.
On the contrary we suppose that $\beta''$ contain $f$ and denote it as $f^{(v'')}$,
where $v''$ is the first coordinate of $f^{(v)}$ in $\beta''$.
If $ 1<v''<2r'+1$,
then $v''=r'+1$ (see Fig. 7 (a)) by Claim 1.
So this copy of $f$ just is the factor $f'$ that contained in $\beta'$.
But it is impossible because its $(n-v''+1)^{th}$
bit has be changed from $a'_{n}$ to $\overline{a'}_{n}$.
Hence $2r'+1< v''<3r'+1$.

So the $(m-v''+1)^{th}$ bit and $(n-v''+1)^{th}$ bit of $f^{(v'')}$ are
$\overline{a'}_{m}$ and $\overline{a'}_{n}$, respectively.
We consider $f^{(u_{n})}$, $f^{(v'')}$ and $f^{(u_{l})}$ (see Fig. 7 (b)).
Let $r_{1}=3r'+1-v''$ and $r_{2}=v''-2r'-1$.
Note that $n-m=r'$ and $r_{1}+r_{2}=r'$.
Equations (1) hold for $t=n-3r'$ and so $\overline{a'}_{n}=a'_{n}$ by Lemma 2.3. It is impossible.
Thus both $\alpha''$ and $\beta''$ contain no $f$ as
factor and so they are 3-critical words of $Q_{d''}(f)$, where $d''=|f|+3r'$.
By $d_{0}=B(f)$ and $\alpha''$ is a factor of $\alpha$, $\alpha''=\alpha$, $\beta''=\beta$ and $d_{0}=d''=|f|+3r'$.

\textbf{Case 2.} $i_{j}\not\in [u_{i_{1}},u_{i_{1}}+|f|-1]$.

If $p=2$, then $\beta$ contains the factor $f^{(u_{i_{1}})}$, a contradiction.
Hence $p\geq3$ and $P=\{1,2,\ldots,p\}\setminus \{1,j\} \neq \emptyset$.
If $i_{k}\not\in [i_{1}+1,u_{i_{1}}+|f|-1]$ for all $k\in P$,
then $\beta$ contains the factor $f^{(u_{i_{1}})}$, a contradiction.
So there must exist $k\in P$ such that $i_{k}\in [i_{1}+1,u_{i_{1}}+|f|-1]$.
Since $j$ is the minimal one such that $i_{1}\in [u_{i_{j}},u_{i_{j}}+|f|-1]$,
$i_{1}\not\in [u_{i_{k}},u_{i_{k}}+|f|-1]$ for any such $k$. So $i_{1}<u_{i_{k}}$.

Let $r=u_{i_{k}}-u_{i_{j}}$ and $s=u_{i_{j}}-u_{i_{1}}$. Obviously $|f|>r+s$.
Let $\alpha'=a'_{1}\cdots a'_{d'}$ be the factor of $\alpha$
that from the $u_{i_{j}}^{th}$ bit to $(u_{i_{k}}+|f|-1)^{th}$ bit of $\alpha$,
where $d'=|f|+r$.
Let $m=i_{k}-u_{i_{j}}+1$ and $n=i_{j}-u_{i_{j}}+1$.
Then $a'_{m}=a_{i_{k}}$ and $a'_{n}=a_{i_{j}}$.
Obviously, $\alpha'$ contains no factor $f$. Denote $\beta'=\alpha'+e_{m}+e_{n}$.
If $\beta'$ does not contain factor $f$, then $\alpha'$ and $\beta'$ is 2-critical words for $Q_{d'}$,
where $d'=|f|+r\leq d-u_{i_{j}}+1<d_{0}$. But this is a contradiction to $d_{0}=B(f)$.
Hence $\beta'$ must contain factor $f$, and so $r\geq2$.
Suppose the first coordinate of this copy of $f$ is $v$, then $1<v\leq r$.
Similar to Claim 1 of  Subcase 1.2 we can show that $r=2r'$ is even,
$n=m+\frac{r}{2}=r'$ and $\overline{a'}_{m}=a'_{n}$.
This indicate that $i_{j}-i_{k}=r'$,
$\overline{a}_{i_{k}}=a_{i_{j}}$ and the factor that start from the $(u_{i_{j}}+r')^{th}$
to the $(u_{i_{j}}+r'+|f|-1)^{th}$ bit of $\alpha+e_{i_{k}}+e_{i_{j}}$ is $f$,
and denote it as $f^{\star}$.
Assume that $i_{1}+r'\neq i_{k}$.
Similar to Claim 4 of Subcase 1.2 we can show that equations (1) hold for
$t=i_{k}-u_{i_{j}}-r'+1$ by considering $f^{\star}$, $f^{(u_{i_{j}})}$ and $f^{(u_{i_{1}})}$, and so $\overline{a}_{i_{k}}=a_{i_{k}}$, a contradiction.
Hence $i_{1}+r'=i_{k}$ and thus we get $\overline{a}_{i_{k}}=a_{i_{1}}$.
Furthermore we consider $f^{(u_{i_{k}})}$, $f^{(u_{i_{j}})}$ and $f^{(u_{i_{1}})}$.
If $s\neq r'$, then equations (1) hold for $t=i_{k}-u_{i_{k}}+1$,
and so $a_{i_{k}}=\overline{a}_{i_{k}}$ by Lemma 2.3, a contradiction.
Hence $s=r'$.

Let $\alpha''=a''_{1}\cdots a''_{d''}$ be the factor of $\alpha$
that start from the $u_{i_{1}}^{th}$ to $(u_{i_{k}}+|f|-1)^{th}$ bit, where $d''=|f|+3r'$.
Let $l'=i_{1}-u_{i_{1}}+1$, $m'=i_{k}-u_{i_{1}}+1$ and $n'=i_{j}-u_{i_{1}}+1$,
then $a''_{l}=a_{i_{1}}$, $a''_{m}=a_{i_{k}}$ and $a''_{n}=a_{i_{j}}$.
It follows that all of $\alpha''+e_{l'}$, $\alpha''+e_{m'}$ and $\alpha''+e_{n'}$ contain factor $f$, and by the above discussions we get $u_{l'}=1$, $u_{n'}=r'+1$ and $u_{m'}=3r'+1$.
Obviously $\alpha''$ does not contain $f$ since it is a factor of $\alpha$.
Let $\beta''=\alpha''+e_{l'}+e_{m'}+e_{n'}$.
Now we show that $\beta''$ contains no factor $f$.
On the contrary we suppose that $\beta''$ contains the factor $f'$
as a copy of $f$ and the first bit of it is $v$.
If $r'<v<3r'+1$, then the only possible case is that $v=2r'+1$,
and so it is just the factor $f$ that contained in $\beta'$ whose the $(2r'+1)^{th}$ bit is $a_{l}$,
but now it is $\overline{a}_{l}$, a contradiction.
Hence $1<v\leq r'$. By considering $f^{(u_{m})}$, $f^{'}$ and $f^{(u_{l})}$,
equations (1) hold for $t=l'-r'$, and so $a'_{l}=\overline{a}_{l}$. It is impossible.

Hence $\beta''$ contains no factor $f$.
Thus $\alpha=\alpha''$,
$\beta=\beta''$ and $\alpha,\beta$ is 3-critical words for $Q_{d_{0}}(f)$,
where $d_{0}=B(f)=|f|+3r'$ and $|f|>3r'$.

Thus by the above discussion we show that if there exist $p$-critical words for $Q_{B(f)}(f)$,
then $p=2$ or $p=3$. The theorem holds.
$\hfill\blacksquare$

By Theorem 3.2, we get $p=2$ or $p=3$.
Obviously, the $p$-critical words $\alpha,\beta$ for $Q_{d_{0}}(f)$ can be determined by the positions of the $p$ copies of $f$ that contained in $\alpha_{j}$, where $j=1,\ldots,p$.
Now we summarize all the possible positions of the copies of $f$ contained in $\alpha_{j}$ for $j=1,\ldots,p$.
It is not difficult to find that if we consider $\alpha^{R}$ and $\beta^{R}$ in Case 2,
then the positions of the copies of $f$ is a special case of Subcase 1.2 such that $i_{3}>|f|$.
So we have the following note.

\trou\noi{\bf Note 3.5}.
The positions of the copies (or their reverse) of $f$ contained in $\alpha_{j}$ only might be the two cases, where $j=1,\ldots,p$ and $p=2$ or $3$.

(1) If $p=2$. Then $u_{i_{1}}=1$ and $u_{i_{2}}=1+r$,
where $1\leq r \leq |f|-2$ and $r+1\leq i_{1}<i_{2}\leq |f|$.
The illustration is shown in Fig. 8 $(a)$.

(2) If $p=3$. Then $u_{i_{2}}=1$, $u_{i_{1}}=2r'+1$, $u_{i_{3}}=3r'+1$,
$i_{2}-i_{1}=r'$, $i_{3}-i_{2}=r'$ and $a_{i_{1}}=\overline{a}_{i_{2}}=a_{i_{3}}$,
where $3r'+1\leq |f|$,
$u_{i_{1}}\leq i_{1}<u_{i_{3}}$ and $u_{i_{3}}\leq i_{2}<i_{3}$.
The illustration is shown in Fig. 8 $(b)$.
Note that it might be $i_{3}\leq |f|$, or $i_{3}>|f|$.

\begin{figure}[h]
\begin{center}
\includegraphics[scale=0.85]{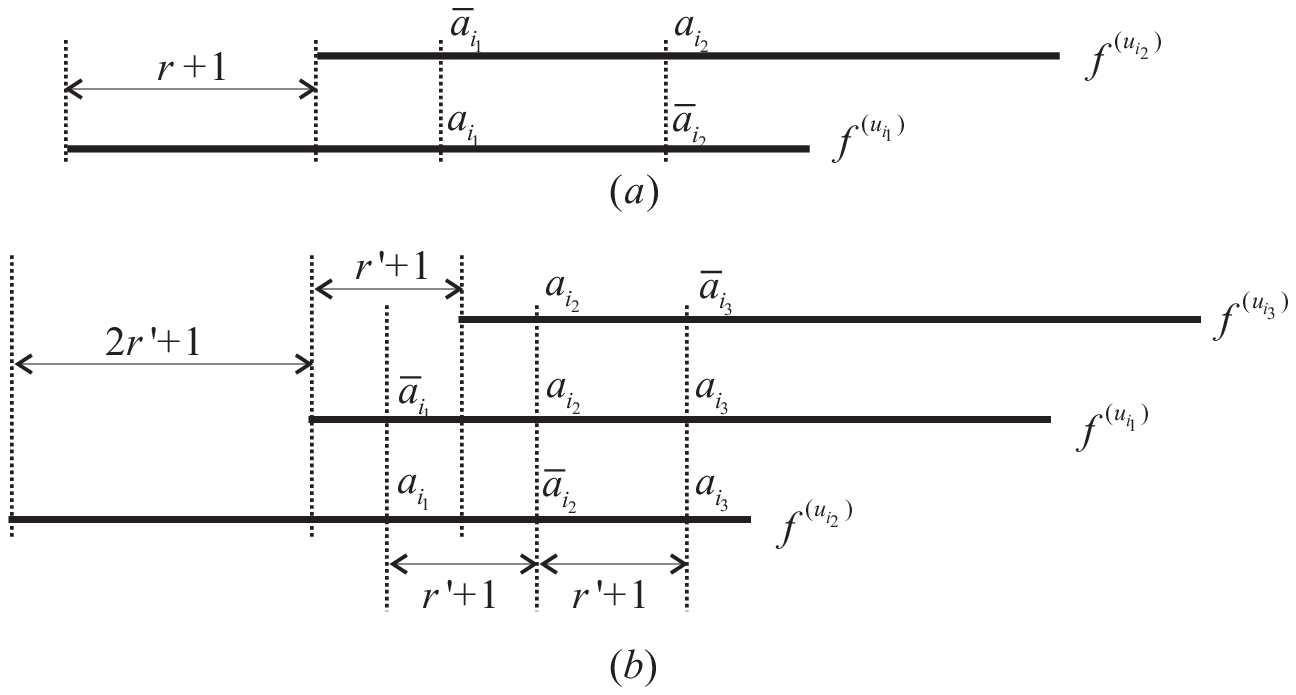}\\
{{\small Fig. 8. Illustration of the positions of $f$ contained in $\alpha_{j}$.}}
\end{center}
\end{figure}

By Note 3.5,
we know that $p$-critical words for $Q_{d_{0}}(f)$ is unique for any bad string $f$.

\trou\noi{\bf Corollary 3.6}. \emph{Let $f$ be a bad word string.
Then there must exist $2$ or $3$-critical words for $Q_{d}(f)$  for any $d\geq B(f)$.}

\trou\noi{\bf Proof}.
Let $\alpha,\beta$ be $p$-critical words for $Q_{B(f)}(f)$.
Then $p=2$ or $p=3$ by Theorem 3.2.
Without loss of generality, suppose the first bit of $f$ is 0.
Then $\alpha'=1^{d-B(f)}\alpha$ and $\beta'=1^{d-B(f)}\beta$ are $p$-critical words for $Q_{d}(f)$ for all $p> B(f)$, where $p=2$ or $p=3$.
$\hfill\blacksquare$

\section{Proofs of Conjectures 1.1 and 1.2}

In this section we give the proofs of Conjectures 1.1 and 1.2.

\trou\noi{\bf Theorem 4.1}. \emph{For any bad string $f$, $B(f)<2|f|$.}

\trou\noi{\bf Proof}. Letting $f$ be any bad string,
there exist $p$-critical words for $Q_{B(f)}(f)$ by Lemma 2.1.
Assume that $\alpha$ and $\beta$ be the $p$-critical words for $Q_{B(f)}(f)$,
then $p=2$ or $p=3$ by Theorem 3.2.
If $p=2$, then $B(f)=|f|+r\leq |f|+|f|-2=2|f|-2$ by Note 3.5 (1).
If $p=3$, then $B(f)=|f|+3r'\leq |f|+|f|-1=2|f|-1$ by Note 3.5 (2).
$\hfill\blacksquare$

To prove Conjecture 1.2, we need the following result.

\trou\noi{\bf Theorem 4.2}.
\emph{If $f$ is a good string, then $ff$ is good.}

\trou\noi{\bf Proof}.
To the contrary we suppose that $ff$ is bad.
Then there exist $\alpha=a_{1}a_{2}\cdots a_{d_{0}}$ and $\beta=b_{1}b_{2}\cdots b_{d_{0}}$ being the $p$-critical words for $Q_{d_{0}}(ff)$ by Lemma 2.1,
where $d_{0}=B(ff)$.
We distinguish two cases since $p=2$ or $p=3$ by Theorem 3.2.

\trou\noi{\bf Case 1}. $p=3$.

Suppose that $a_{i}\neq b_{i}$, $a_{j}\neq b_{j}$ and $a_{k}\neq b_{k}$, where $i<j<k$.
Without loss of generality,
set $\alpha_{1}=\alpha+e_{i}$, $\alpha_{2}=\alpha+e_{j}$ and $\alpha_{3}=\alpha+e_{k}$ $\not\in V(Q_{d_{0}}(ff))$.
It follows that $ff$ is a factor of $\alpha_{t}$, $t=1,2,3$.
Denote the copy of $ff$ contained in $\alpha_{1}$, $\alpha_{2}$ and $\alpha_{3}$ with $f^{(2)}f^{(2')}$, $f^{(3)}f^{(3')}$ and $f^{(1)}f^{(1')}$, respectively, where $f^{(t)}=f^{(t')}=f$, $t=1,2,3$.
By Note 3.5 (2) we suppose the first coordinates of $f^{(1)}f^{(1')}$, $f^{(2)}f^{(2')}$
and $f^{(3)}f^{(3')}$ are $ u_{k}=3r+1, u_{i}=2r+1,u_{j}=1$, respectively.
Furthermore $3r+1\leq |ff|=2|f|$, $j-i=k-j=r$ and $a_{i}=\overline{a}_{j}=a_{k}$.
Clearly, $r\neq |f|$ and $r<2|f|-1$.
If $|f|+1\leq r<2|f|-1$, then $2r\geq 2|f|+2>2|f|=|ff|$, a contradiction to that $3r+1\leq |ff|$.
Thus $1\leq r<|f|$.
By Note 3.5 (2), $3r+1<k\leq 2r+2|f|$.
Now we distinguish three subcases: $3r+1< k\leq 2r+|f|$, $2r+|f|+1\leq k\leq 3r+|f|$
 and $3r+|f|+1\leq k\leq 2r+2|f|$ .

\trou\noi{\bf Subcase 1.1}. $3r+1< k\leq 2r+|f|$.

This case is shown in  Fig. 9 (a).
We consider the factors $f^{(2')}$, $ f^{(1)}$ and $f^{(2)}$.
Their first coordinates are $2r+|f|+1$,
 $3r+1$ and $2r+1$ (see Fig. 9 (b)), respectively, and they satisfy:

$ f^{(2')}_{t}= f^{(1)}_{t+|f|-r}$ for all $t\in \{1,\ldots,r\}$;

$ f^{(1)}_{t}= f^{(2)}_{t+r}$ for all $t\in \{1,\ldots,|f|-r\} \setminus \{k-3r\}$;

$ f^{(1)}_{k-3r}=\overline{a}_{k}$ and $ f^{(2)}_{k-2r}= a _{k}$.

\begin{figure}[h]
\begin{center}
\includegraphics[scale=0.85]{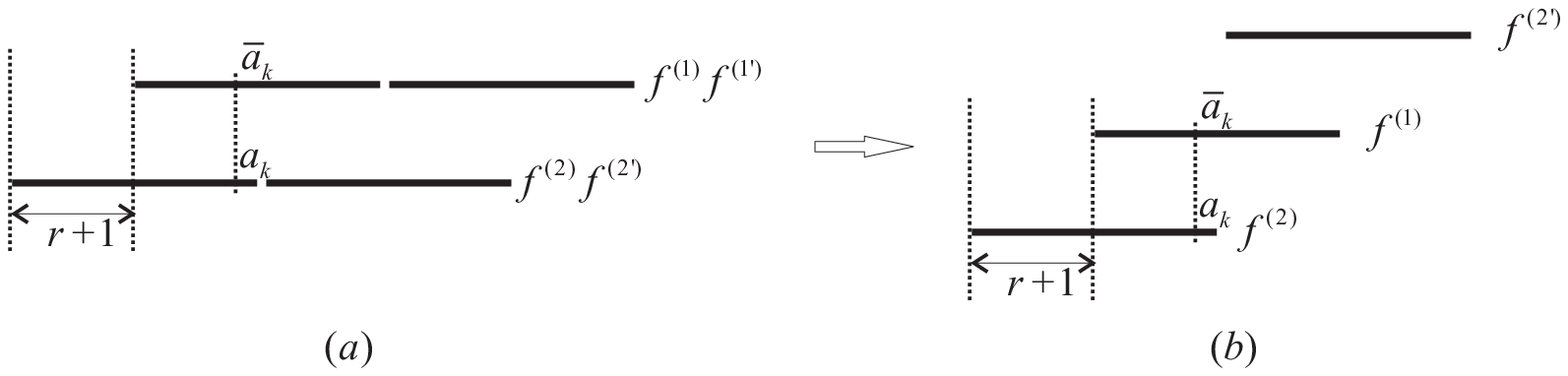}\\
{{\small Fig. 9. Illustration of the proof of Subcase 1.1 of Theorem 4.2}}
\end{center}
\end{figure}

Suppose that $g$ is the greatest common divisor of $r$ and $|f|-r$,
and $t_{0}$ is the integer such that $1\leq t_{0}\leq g$,
$t_{0}\leq k-3r$ and $k-t_{0}\equiv 0$ (mod~$g$).
Then by Corollary 2.5, the equations (2) and (3) hold for $t=t_{0}$,
and so $\overline{a}_{k}=f^{(1)}_{k-3r}=f^{(2)}_{k-2r}=a _{k}$,
a contradiction.

\trou\noi{\bf Subcase 1.2}. $2r+|f|+1\leq k\leq 3r+|f|$.

This case is shown in Fig. 10 (a).
Now we consider $f^{(2')}$, $ f^{(1)}$ and $f^{(2)}$ as shown in Fig. 10 (b).
Their first coordinates are $2r+|f|+1$, $3r+1$ and $2r+1$, respectively,
and they satisfy:

$ f^{(2')}_{t}= f^{(1)}_{t+|f|-r}$ for all $t\in \{1,\ldots,r\}\setminus\{k-2r-|f|\}$,

$ f^{(2')}_{k-2r-|f|}= a_{k}$, $ f^{(1)}_{k-r}=\overline{a} _{k}$ and

$ f^{(1)}_{t}= f^{(2)}_{t+r}$ for all $t\in \{1,\ldots,|f|-r\}$.

\begin{figure}[h]
\begin{center}
\includegraphics[scale=0.85]{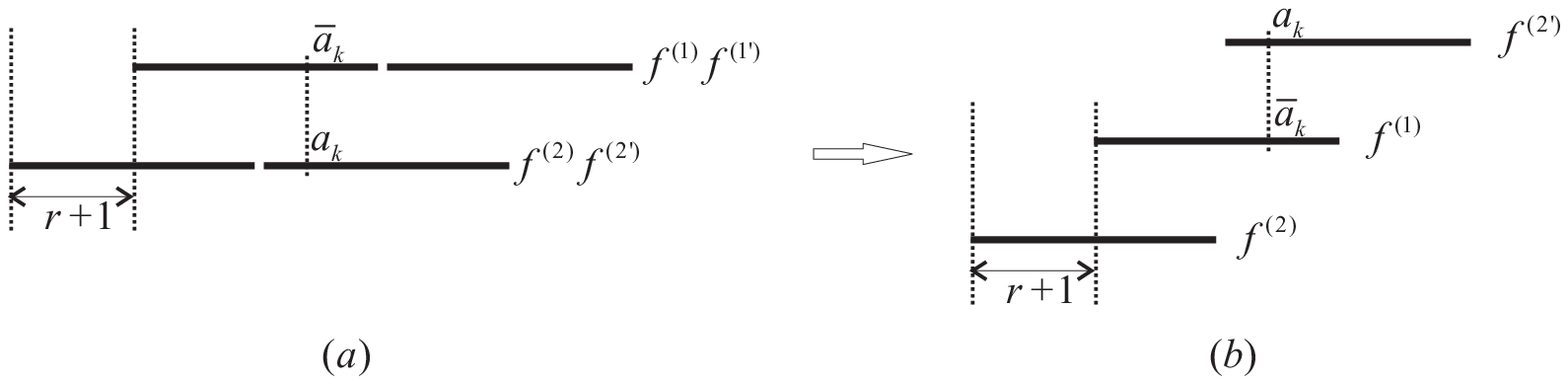}\\
{{\small Fig. 10. Illustration of the proof of Subcase 1.2 of Theorem 4.2.}}
\end{center}
\end{figure}

Let $t_{0}$ be the positive integer such that $1\leq t_{0}\leq g$, $t_{0}\leq k-2r-|f|$ and $k-t_{0}\equiv 0$ (mod $g$), where $g$ is the greatest common divisor of $r$ and $|f|-r$.
By Corollary 2.5, the equations (2) and (3) hold for $t=t_{0}$,
and so $ \overline{a}_{k}= f^{(2')}_{k-2r-|f|}= f^{(1)}_{k-r}=a _{k}$, contradiction.

\trou\noi{\bf Subcase 1.3}. $3r+|f|+1\leq k\leq 2r+2|f|$.

This case is shown in Fig. 11 (a).
Now we consider $f^{(1')}$, $ f^{(2')}$ and $f^{(1)}$ as shown in Fig. 11 (b).
Their first bits start from $3r+|f|+1$,~$2r+|f|+1$ and $3r+1$, respectively,
and they satisfy:

$ f^{(1')}_{t}= f^{(2')}_{t+r}$,  for all $t\in \{1,\ldots,|f|-r\} \setminus\{k-3r-|f|\}$,

$ f^{(1')}_{k-3r-|f|}=\overline{a}_{k}$, $ f^{(2')}_{k-2r-|f|}= a _{k}$ and

$ f^{(2')}_{t}= f^{1'}_{t+|f|-r}$ for all $t\in \{1,\ldots,r\}$.

\begin{figure}[h]
\begin{center}
\includegraphics[scale=0.85]{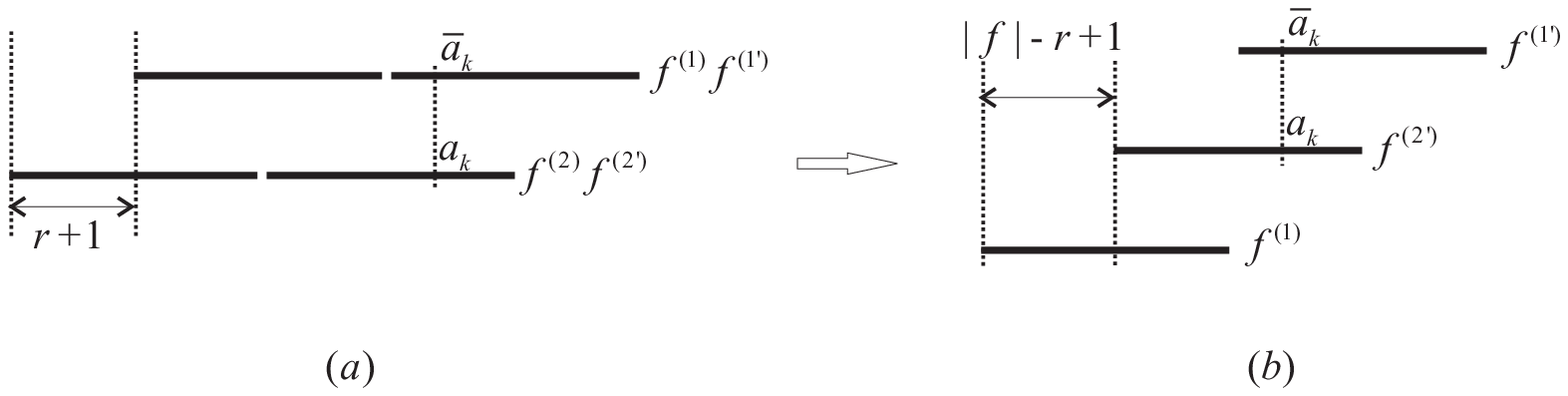}\\
{{\small Fig. 11. Illustration of the proof of Subcase 1.3 of Theorem 4.2}}
\end{center}
\end{figure}

Let $t_{0}$ be the positive integer such that $1\leq t_{0}\leq g$, $t_{0}\leq k-3r-|f|$
and $k-t_{0}\equiv 0$ (mod~$g$), where $g$ is the greatest common divisor of $r$ and $|f|-r$.
By Corollary 2.5, the equations (2) and (3) hold for $t=t_{0}$,
and so $\overline{a}_{k}=f^{(1')}_{k-3r-|f|}=f^{(2')}_{k-2r-|f|}=a _{k}$, a contradiction.

Thus we get that if $ff$ is a bad string, then $Q_{d_{0}}(ff)$ has no $3$-critical words.

\trou\noi{\bf Case 2}. $p=2$.

Suppose that $a_{i}\neq b_{i}$ and $a_{j}\neq b_{j}$, where $i<j$.
Then $\alpha_{1}=\alpha+e_{i}$ and $\alpha_{2}=\alpha+e_{j}$ $\not\in V(Q_{d_{0}}(ff))$.
Denote the factor $ff$ that contained in $\alpha_{1}$ and $\alpha_{2}$ with $f^{(1)}f^{(1')}$ and $f^{(2)}f^{(2')}$, respectively, where $f^{(t)}=f^{(t')}=f$ for $t=1,2$.
By Note 3.5 (1), we may suppose $u_{j}=1$ (the case $u_{i}=1$ can be treated analogously).
Letting $r=u_{i}-u_{j}$, $1\leq r\neq |f|$ and $r<2|f|-1$.
By $1\leq r<|f|$ or $|f|+1\leq r<2|f|-1$, we distinguish two subcases.

\trou\noi{\bf Subcase 2.1}. $1\leq r<|f|$.

Clearly, $r+1\leq i<j\leq2|f|$. We consider the possible positions of $i$ and $j$ further.

Suppose $r+1\leq i\leq|f|+r$ and $|f|+r+1\leq j\leq 2|f|$.
Then we consider the factors $f^{(2')}$, $f^{(1)}$ and $f^{(2)}$.
Their first bits start from $|f|+1$, $r+1$ and $1$, respectively.
By similar discussion as in subcase 1.1, $\overline{a}_{i}=a_{i}$, obviously a contradiction.
Thus it does not hold that $r+1\leq i\leq|f|+r$ and $|f|+r+1\leq j\leq 2|f|$.

Assume that $r+1\leq i\leq|f|$ and $|f|+1\leq j\leq2|f|$.
We consider the factors $f^{(1')}$, $f^{(2)}$ and $f^{(1)}$.
By a similar discussion as in Subcase 1.2,
we get $\overline{a}_{j}=a_{j}$, a contradiction.
Thus $r+1\leq i\leq|f|$ and $|f|+1\leq j\leq2|f|$ does not hold.

So the possible positions of $i$ and $j$ only might be:
$r+1\leq i< j\leq|f|$,
$|f|+1\leq i<j\leq|f|+r$,
or $|f|+r+1\leq i<j\leq 2|f|$.

First we consider the case $r+1\leq i< j\leq|f|$.
If $j-i\not\equiv 0$ (mod $g$),
then $a_{i}=\overline{a}_{i}$ by considering $f^{(2')}$, $f^{(1)}$ and $f^{(2)}$, a contradiction.
Hence $j-i \equiv 0$ (mod $g$). By Lemma 2.6 $(\romannumeral2)$, $a_{i}=a_{j}$ .
Let $\alpha'=a_{1}\cdots a_{|f|+r},\beta'=b_{1}\cdots b_{|f|+r}$ be the strings obtained from $\alpha,\beta$ by deleting the factor that start from the $(|f|+r+1)^{th}$ bit to the last bit, respectively.
It follows that both $\alpha'+e_{i}$ and $\alpha'+e_{j}$ contain $f^{(1)}$ and $f^{(2)}$ as factor,  and their first bits are $ u_{i}=r+1$ and $ u_{j}=1$, respectively.
Obviously if $r=1$, then both $\alpha'$ and $\beta'$ contain no factor $f$.
Assume that $r\geq 2$ and $\alpha'$ contains factor $f'$ as a copy of $f'$.
Clearly the first bit $v$ of $f'$ satisfy $1<v\leq r$.
Now we consider $f^{(1)},f'$ and $f^{(2)}$.
It is easy to see that equations (1) hold for $t= j-r$,
and so $a_{j}=\overline{a}_{j}$, a contradiction.
Assume that $\beta'$ contain $f$ as factor, and denote it with $f^{''}$,
Clearly the first bit $v'$ of $f''$ satisfy $1<v'\leq r$.
Now we consider $f^{(1)},f^{''}$ and $f^{(2)}$.
By $a_{i}=a_{j}$ and equations (1) hold for $t=j-r$,
we get $a_{j}=\overline{a}_{j}$, a contradiction.
Thus both $\alpha'$ and $\beta'$ contain no factor $f$.
So $\alpha'$ and $\beta'$ are 2-critical words for $Q_{d'}(f)$, where $d'=|f|+r$.
Hence $f$ is a bad string, a contradiction since we assume that $f$ is a good string.

Now we consider the case $|f|+1\leq i<j\leq|f|+r$.
We consider $f^{(2')}$, $f^{(1)}$ and $f^{(2)}$.
If $i-j\not\equiv 0\pmod g$, then $a_{i}=\overline{a}_{i}$ by Corollary 2.5, a contradiction.
Hence $i-j \equiv 0\pmod g$, and so by Lemma 2.6 $(\romannumeral3)$, $a_{i}= a_{j}$.
Let $\alpha' $ and $\beta' $ be the factors of $\alpha$ and $\beta$ that start from the $(r+1)^{th}$ bit to the $(2|f|)^{th}$ bit, respectively.
Now we claim that both $\alpha' $ and $\beta'$ contain no $f$ as factor.
Obviously it holds for $r=1$.
By contrary assume that $\alpha'$ contain factor $f$ for $r\geq 2$,
and denote this factor with $f'$ and its first coordinate is $v$.
Hence by considering $f^{(2')}$, $f'$, $f^{(1)}$ and $a_{i}= a_{j}$,
equations (1) hold for $t=i-|f|+r$,
and so $a_{i}= \overline{a}_{i}$, a contradiction.
Assume that $\beta'$ contain factor $f$ for $r\geq 2$,
and denote this factor with $f''$ with first coordinate $v'$.
Hence by considering $f^{(2')}$, $f''$ and $f^{(1)}$,
equations (1) hold for $t=i-|f|+r$,
and so $a_{i}= \overline{a}_{i}$, a contradiction.
Thus $\alpha'$ and $\beta'$ are 2-critical words for $Q_{d'}(f)$, where $d'=2|f|-r$.
So $f$ is a bad string. This is a contradiction since we assume that $f$ is a good string.

Finally we consider the case $|f|+r+1\leq i<j\leq 2|f|$.
We consider $f^{(1')}$, $f^{(2')}$ and $f^{(1)}$.
If $i-j\not\equiv 0\pmod g$,
then we get $a_{i}=\overline{a}_{i}$ by Corollary 2.5, a contradiction.
Hence $i-j \equiv 0\pmod g$. By Lemma 2.6 $(\romannumeral3)$, we get $a_{i}= a_{j}$ .
Let $\alpha' $ and $\beta' $ be the factors of $\alpha$ and $\beta$ that start from
the $(|f|+1)^{th}$ bit to the last bit, respectively.
We claim that both $\alpha'$ and $\beta'$ contain no factor $f$.
Obviously it holds for $r=1$.
By contrary we suppose $\alpha'$ contain factor $f$ and denote it as $f'$.
Let the first coordinate of $f'$ be $v$.
We consider $f^{(1')}$, $f'$ and $f^{(2')}$.
Obviously equations (1) hold for $t=i-r$, and so $a_{i}= \overline{a}_{i}$, a contradiction.
By contrary we suppose $\beta'$ contain factor $f$ and denote it as $f''$.
Let the first coordinate of $f''$ be $v'$.
We consider $f^{(1')}$, $f''$ and $f^{(2')}$.
By $a_{i}= a_{j}$ and equations (1) holding for $t=i-r$,
and so $a_{i}= \overline{a}_{i}$, a contradiction.
Thus $\alpha$ and $\beta$ are 2-critical words for $Q_{d'}(f)$, where $d'=|f|+r$.
So $f$ is a bad string. This is a contradiction since we assume that $f$ is a good string.

Hence in the subcase we get that if $f$ is good, then $ff$ is good.

\trou\noi{\bf Subcase 2.2}. $|f|+1\leq r<2|f|-1$.

Obviously, $r+1\leq i<j\leq 2|f|$.
Let $\alpha'=a'_{1}\cdots a'_{d'}$ and $\beta'=b'_{1}\cdots b'_{d'}$
be the factors of $\alpha$ and $\beta$ that start from the $(|f|+1)^{th}$
 bit to the $(|f|+r)^{th}$, respectively, where $d'=r$.
Letting $ l =i-|f|$ and $ m = j-|f|$, $a'_{l}=a_{i}$ and $a'_{m}=a_{j}$.
Obviously both $\alpha'+e_{i}$ and $ \alpha'+e_{j}$ contain $f$ as factor.
Letting $r'=r-|f|$,
the first bits of the two copies $f^{(u_{l})},f^{(u_{m})}$ of
$f$ are $u_{m}=1, u_{l}=r'+1$, respectively.

Obviously, if $r'=1$, then both $\alpha'$ and $\beta'$ contain no $f$ as factor.
Now assume that $r'\geq 2$.
Suppose that $\alpha' $ contains $f$ as a factor and denote it as $f^{'}$.
Then the first bit $v$ of $f^{'}$ satisfies $1<v\leq r'$.
By consider $f^{(1)}$, $f^{'} $ and $f^{(r'+1)}$ equations (1) hold for $t=l-r'$,
and so $\overline{a'}_{l}=a'_{l}$, a contradiction.
Hence $\alpha' $ contains no $f$ as factor.

Assume that $\beta' $ contains factor $f$ and denote it as $f^{''}$.
Then the first bit $v'$ of $f^{''}$ satisfies $1<v'\leq r'$.
By Claim 1 of Subcase 1.2 of Theorem 3.2,
$r'$ must be even, $\frac{r'}{2}+l=m$ and $\overline{a'}_{l}=a'_{m}$.
Letting $r'=2s$, $v'=\frac{r'}{2}+1=s+1$.

Now by considering $f^{(1)}$, $f^{''}$ and $f^{(r'+1)}$,
we characterize the structure of $f$, $\alpha'$ and $\beta'$.
Letting $k_{1}$ be the positive integer such that $l-k_{1}s\geq1$ and $l-(k_{1}+1)s<1$, $k_{1}\geq2$.
Letting $t_{1}=l-k_{1}s-1$, $0\leq t_{1}\leq s-1$.
Letting $k_{2}$ be the positive integer such that $m+k_{2}s\leq |f|+2s$ and $m+(k_{2}+1)s>|f|+2s$, $k_{2}\geq2$.
Letting $t_{2}=|f|+2s-m-k_{1}s$, $0\leq t_{2}\leq s-1$.
Let $\mu$ be the factor of $\alpha'$ that start from the $(l+1)^{th}$ bit to the $(m-1)^{th}$ bit.
Clearly, $|\mu|=s-1$.
Let $ \rho$ be the factor of $\alpha'$ that start from the first bit to the $t_{1}^{th}$ bit,
and $\sigma$ be the factor of the last $t_{2}$ bits.
Note that $\mu$, $\rho$ and $\sigma$ might be null strings.
Letting $ a'_{l}=x$ and $ a'_{m}=\overline{x}$,
the $(l-ts)^{th}$ bit of $\alpha'$ is $\overline{x}$, where $t=1,\ldots,k_{1}$,
the $(m+ts)^{th}$ bit of $\alpha'$ is $x$, where $t=1,\ldots,k_{2}$,
and for the other bits $k',k''$ such that $k'\equiv k''$ and $k'\not\equiv l\pmod {r'}$,
 $a_{k'}=a_{k''}$.
Hence $\sigma$ and $\rho^{R}$ are factors of $\mu$ and $ \mu^{R}$ that start from the first bit, respectively.
Note that $\beta'=\alpha' +e_{l} +e_{m}$,
and $f$ can be obtained by deleting the first $s$ bits and the last $s$ bits of $\beta'$.
We get

$\alpha'=\rho\overbrace{\overline{x}\mu\cdots\overline{x}\mu}^{k_{1}}x\mu\overline{x}
\overbrace{\mu x\cdots \mu x}^{k_{2}}\sigma$,

$\beta'=\rho\overbrace{\overline{x}\mu\cdots\overline{x}\mu}^{k_{1}}\overline{x}\mu x
\overbrace{\mu x\cdots \mu x}^{k_{2}}\sigma$ and

$f=\rho\overbrace{\overline{x}\mu\cdots\overline{x}\mu}^{k_{1}-1}\overline{x}\mu x
\overbrace{\mu x\cdots\mu x}^{k_{2}-1}\sigma$,

Now we construct strings $\alpha''$ and $\beta''$ as following,
where $\alpha''$ and $\beta''$ are obtained by add $\overbrace{x\mu\cdots x\mu}^{k_{2}-2}$
behind the $l^{th}$ bit and
 $\overbrace {\mu \overline{x}\cdots \mu \overline{x}}^{k_{1}-2}\mu x$ behind the $m^{th}$ bit of $\alpha'$ and $\beta'$, respectively.

~$\alpha''=\rho\overbrace{\overline{x}\mu\cdots\overline{x}\mu}^{k_{1}}
\overbrace {x\mu\cdots x\mu}^{k_{2}-2}
x\mu\overline{x}
\overbrace {\mu \overline{x}\cdots \mu \overline{x}}^{k_{1}-2}\mu x
\overbrace{\mu x\mu\cdots x}^{k_{2}}\sigma$,

~$\beta''=\rho\overbrace{\overline{x}\mu\cdots\overline{x}\mu}^{k_{1}}
 \overbrace {x\mu\cdots x\mu}^{k_{2}-2}
\overline{x}\mu x
 \overbrace {\mu \overline{x}\cdots \mu \overline{x}}^{k_{1}-2}\mu \overline{x}
\overbrace{\mu x\mu\cdots x}^{k_{2}}\sigma$.

Obviously, $|\alpha''|=|\beta''|=d''=|f|+(k_{1}+k_{2}-1)s$.
Note that $ \alpha''$ and $\beta'' $ differ in exactly in the $l_{0}^{th}$,
$m_{0}^{th}$ and $n_{0}^{th}$ bits,
where $l_{0}=l+(k_{2}-2)s$,
$m_{0}=l+(k_{2}-1)s$ and $n_{0}=l+(k_{2}+k_{1}-3)s$.
Hence $H(\alpha'', \beta'')=3$.
Obviously all of $\alpha''+e_{l_{0}}$, $\alpha''+e_{m_{0}}$ and $\alpha''+e_{n_{0}}$
contain factor $f$, and both $\alpha''$ and $\beta''$ contain no factor $f$.
Thus $\alpha''$ and $\beta''$ are 3-critical words for $Q_{d''}(f)$.
So $f$ is a bad string, a contradiction since $f$ is good.

Hence by the above discussion we know that if $ff$ is bad, then $f$ is bad.
This completes the proof.
$\hfill\blacksquare$

\trou\noi{\bf Theorem 4.3}. \emph{If $Q_{d}(f)$ is an isometric subgraph of $Q_{d}$,
then $Q_{d}(ff)$ is an isometric subgraph of $Q_{d}$.}

\trou\noi{\bf Proof}. 
First we suppose that $f$ is a bad string.
If $Q_{d}(f)\hookrightarrow Q_{d}$, then $d\leq 2|f|-1$ by Theorem 4.1.
Hence $|ff|=2|f|>d$.
Clearly, $Q_{d}(ff)\cong Q_{d}$ and so $Q_{d}(ff)\hookrightarrow Q_{d}$.

Now we suppose $ f $ is good. Then $ ff $ is good by Theorem 4.2.
So $Q_{d}(ff)\hookrightarrow Q_{d}$.
$\hfill\blacksquare$




\begin{thebibliography}{6}
\bibitem{dt} E. Ded\'{o}, D. Torri, N. Z. Salvi,
The observability of the Fibonacci and the Lucas cubes,
Discrete Math. 255 (2002) 55--63.
\bibitem{Gre} P. Gregor,
Recursive fault-tolerance of Fibonacci cube in hypercubes,
Discrete Math. 306 (2006) 1327--1341.
\bibitem{whs} W.-J. Hsu,
Fibonacci cubes - a new interconnection topology,
IEEE Trans. Parallel Distrib. Syst. 4 (1) (1993) 3--12.
\bibitem{asy1} A. Ili\'{c}, S. Klav\v{z}ar, Y. Rho,
Generalized Fibonacci cubes,
Discrete Math. 312 (2012) 2--11.
\bibitem{asy} A. Ili\'{c}, S. Klav\v{z}ar, Y. Rho, The index of a binary word,
Theoret. Comput. Sci. 452 (2012) 100--106.
\bibitem{ik} W. Imrich, S. Klav\v{z}ar,
Product graphs: structure and recognition,
Wiley, New York, 2000.
\bibitem{sk1} S. Klav\v{z}ar,
On median nature and enumerative properties of Fibonacci-like cubes,
Discrete Math. 299 (2005) 145--153.
\bibitem{sk} S. Klav\v{z}ar,
Structure of Fibonacci cubes: a survey,
J. Comb. Optim. 25 (2013) 505--522.
\bibitem{sks} S. Klav\v{z}ar, S. Shpectorov,
Asymptotic number of isometric generalized Fibonacci cubes,
Europ. J. Combin. 33 (2012) 220--226.
\bibitem{sp} S. Klav\v{z}ar, P. \v{Z}igert,
Fibonacci cubes are the resonance graphs of fibonaccenes,
Fibonacci Quart. 43 (3) (2005) 269--276.
\bibitem{liu} J. Liu, W.-J. Hsu, M. J. Chung,
Generalized Fibonacci cubes are mostly Hamiltonian,
J. Graph Theory 18 (8) (1994) 817--829.
 \bibitem{ms} E. Munarini, N. Z. Salvi,
 Structural and enumerative properties of the Fibonacci cubes,
Discrete Math. 255 (2002) 317--324.
\bibitem{nzs} N. Z. Salvi,
On the existence of cycles of every even length on generalized Fibonacci cubes,
Matematiche (Catania) 51 (1996) 241--251.
\bibitem{tara} A. Taranenko, A. Vesel,
Fast recognition of Fibonacci cubes,
Algorithmica 49 (2) (2007) 81--93.
\bibitem{zoy} H. Zhang, L. Ou, H. Yao,
Fibonacci like cubes as Z-transformation graphs,
Discrete Math. 309 (2009) 1284--1293.
\end{thebibliography}
\end{document}